\documentclass[11pt]{amsart}

\usepackage{amssymb,amsthm,amsmath}
\usepackage[numbers,sort&compress]{natbib}
\usepackage[margin=1in]{geometry}
\usepackage{color}
\usepackage{graphicx}
\usepackage{tikz}
\usepackage[colorlinks,
linkcolor=red,
anchorcolor=blue,
citecolor=green
]{hyperref}

\def\e{\varepsilon}

\let\newpf\proof \let\proof\relax 
\newenvironment{pf}{\newpf[\proofname]}{\qed\endtrivlist}

\newcommand{\ba}{\overline{A}}

\def\be{\begin{equation}}
\def\ee{\end{equation}}

\def\ba{{\begin{align}}}
\def\ea{{\end{align}}}

\def\bm{\begin{matrix}}
\def\em{\end{matrix}}

\def\0{{\mathbf 0}}

\newtheorem{Theorem}{Theorem}[section]
\newtheorem{Lemma}{Lemma}[section]
\newtheorem{Proposition}{Proposition}[section]
\newtheorem{Corollary}{Corollary}[section]
\theoremstyle{definition}
\newtheorem{Remark}{Remark}[section]

\newtheorem{Definition}{Definition}[section]

\numberwithin{equation}{section}

\theoremstyle{definition}

\newtheorem{definition}{Definition}[section]

\newcommand{\C}{{\mathbb C}}

\newcommand{\Q}{{\mathbb Q}}
\newcommand{\R}{{\mathbb R}}
\newcommand{\T}{{\mathbb T}}

\newcommand{\Z}{{\mathbb Z}}

\def\B0{{\bold{0}}}

\catcode`\@=12

\def\Empty{}
\newcommand\oplabel[1]{
  \def\OpArg{#1} \ifx \OpArg\Empty {} \else
    \label{#1}
  \fi}

\newcommand{\comm}[1]{}
\newcommand{\comment}[1]{}

\begin{document}
\title[]{ H\"older regularity of the integrated density of states for quasi-periodic long-range operators on $\ell^2(\Z^d)$}
\author{Lingrui Ge}
\address{Department of Mathematics, University of California Irvine, CA, 92697-3875, USA}
\email{lingruig@uci.edu}

\author {Jiangong You}
\address{
Chern Institute of Mathematics and LPMC, Nankai University, Tianjin 300071, China} \email{jyou@nankai.edu.cn}

\author{Xin Zhao}
\address{
Department of Mathematics, Nanjing University, Nanjing 210093, China
and
Department of Mathematics, University of California Irvine, CA, 92697-3875, USA}
 \email{njuzhaox@126.com}

%\author {Jiangong You}
%\institute{
%Chern Institute of Mathematics and LPMC, Nankai University, Tianjin 300071, China} 

%\author{Xin Zhao}
%\address{
%Department of Mathematics, Nanjing University, Nanjing 210093, China
%\and
%Department of Mathematics, University of California Irvine, CA, 92697-3875, USA}
 %\email{njuzhaox@126.com}

\begin{abstract}
We prove  the  H\"older continuity of the integrated density of states for a class of  quasi-periodic long-range operators on $\ell^2(\Z^d)$ with large trigonometric polynomial potentials and Diophantine frequencies. Moreover, we give  the H\"older exponent in terms of the cardinality of the level sets  of the potentials,  which improves, in the perturbative regime, the result obtained by Goldstein and Schlag \cite{gs2}. Our approach is a combination of Aubry duality, generalized Thouless formula and the regularity of the Lyapunov exponents of analytic quasi-periodic $GL(m,\C)$ cocycles which is proved by quantitative almost reducibility method. 
\end{abstract}

\maketitle
\section{Introduction}
In this paper, we consider the following quasi-periodic long-range operator on $\ell^2(\Z^d)$ with trigonometric polynomial potentials:
\begin{equation}\label{1.2}
(L^{\lambda W}_{V,\alpha,\theta}u)_n=\sum\limits_{k\in\Z^d}\widehat{V}_ku_{n-k}+\lambda W(\theta+\langle n,\alpha\rangle)u_{n}, \ \ n\in \Z^d,
\end{equation}
where  $W(\theta)=\sum_{k=-m}^m\widehat{W}(k)e^{2\pi ik\theta}$ is a real trigonometric polynomial and  $V(x)=\sum_{k\in\Z^d}\widehat{V}(k)e^{2\pi i\langle k,x\rangle}$ is a real analytic function. $W$, $\theta\in \T$,  $\alpha\in\R^d$ are called  the potential,  the phase and  the frequency respectively. We always assume that $\{1,\alpha_1,\cdots,\alpha_d\}$ are independent over $\Q$. \eqref{1.2} includes and relates to several interesting quasi-periodic models.\\
{\bf Example 1.} Taking $V(x)=\sum\limits_{i=1}^d2\cos 2\pi x_i$, it is the quasi-periodic Schr\"odinger operator on $\ell^2(\Z^d)$,
\begin{equation}\label{highschrodinger-1}
H= \Delta +  \lambda W(\theta+\langle n,\alpha\rangle)\delta_{nn'},
\end{equation}
where $\Delta$ is the usual Laplacian on $\Z^d$ lattice.\\
{\bf Example 2.} Taking $d=1$, $V(x)=2\cos 2\pi x$, it is the one-frequency quasi-periodic Schr\"odinger operator on $\ell^2(\Z)$,   
\begin{equation}\label{schrodinger-1}
(H_{\lambda W,\alpha,\theta}u)(n)= u_{n+1}+u_{n-1} + \lambda W(\theta+n\alpha)u(n).
\end{equation}
{\bf Example 3.} Taking $d=m=1$, $V(x)=W(x)=2\cos 2\pi x$, it is the famous almost Mathieu operator,
\begin{equation}\label{amo}
(H_{\lambda,\alpha,\theta}u)(n)= u_{n+1}+u_{n-1} +2\lambda\cos 2\pi(\theta+n\alpha)u(n).
\end{equation}
{\bf Example 4.} The Aubry dual of \eqref{1.2} is the following multi-frequency quasi-periodic finite-range operator on $\ell^2(\Z)$,
\begin{equation}\label{schrodinger-2}
(L^{\lambda^{-1} V} _{W,\alpha,x}u)(n)= \sum\limits_{k=-m}^m\widehat{W}_k u_{n-k} +\lambda^{-1} V(x+n\alpha)u(n).
\end{equation} 
{\bf Example 5.} The special case of (\ref{schrodinger-2}) is the following quasi-periodic Schr\"odinger operator with $d$ frequencies,
\begin{equation}\label{schrodinger-3}
(H _{\lambda^{-1}V,\alpha,x}u)(n)=  u_{n+1}+u_{n-1} +\lambda^{-1} V(x+n\alpha)u(n).
\end{equation} 
The Integrated Density of States (IDS) is a quantity of fundamental importance for models in condensed matter physics. It is defined in a uniform way for all quasi-periodic family self-adjoint operators $(L_x)_{x\in\T^d}$ by
$$
N(E)=\int_{\T^d}\mu_{x}(-\infty,E]dx,
$$
where $\mu_x$ is the spectral measure associated with $L_x$ and $\delta_0$. Roughly speaking, the density of states measure $N([E_1,E_2])$ gives the ``number of states per unit volume" with energy between $E_1$
and $E_2$.

The regularity of IDS is a fascinating subject in the spectral theory of quasi-periodic operators, especially for the absolute continuity  \cite{ad,avila1,aj} and the H\"older regularity \cite{amor,aj,gs1,gs2}. The regularity of IDS is also closely related to many other topics in the spectral theory of quasi-periodic operators. For example, the absolute continuity of IDS is closely related to purely absolutely continuous spectrum in the regime of zero Lyapunov exponent \footnote{See section 2.1 for the definition of Lyapunov exponent of Schr\"odinger cocycles.} \cite{kotani,damanik}. H\"older continuity of IDS is closely related to homogeneity of the spectrum \cite{dgl,dgsv,lyzz}.

\subsection{Regularity of IDS}
We first review the previous results on H\"older regularity of IDS for one dimensional quasi-periodic Schr\"odinger operators. For (\ref{schrodinger-3})  with real analytic potential $V$ (equivalently for its dual operator (\ref{1.2}) with $W=2\cos2\pi\theta$), $\frac{1}{2}$-H\"older continuity of IDS was obtained by Amor \cite{amor} for sufficiently large $\lambda$  and Diophantine $\alpha$ \footnote{  Recall that $\alpha \in\R^d$ is called {\it Diophantine} if there are $\kappa>0$ and $\tau>d-1$ such that $\alpha \in {\rm DC}_d(\kappa,\tau)$, where
\begin{equation}\label{dio}
{\rm DC}_d(\kappa,\tau):=\left\{\alpha \in\R^d:  \inf_{j \in \Z}\left|\langle n,\alpha  \rangle - j \right|
> \frac{\kappa}{|n|^{\tau}},\quad \forall \, n\in\Z^d\backslash\{0\} \right\},
\end{equation}
and ${\rm DC}_d:=\bigcup_{\kappa>0,\tau>d-1} {\rm DC}_d(\kappa,\tau)$ is a full Lebesgue measure set.}.  For the one-frequency case, i.e., operator \eqref{schrodinger-1} with $W$ being a real analytic potential, Amor's result was extended by Avila-Jitomirskaya \cite{aj} to the non-perturbative regime (actually to the almost reducible regime). The  H\"{o}lder continuity of IDS in the positive Lyapunov exponent regime was proved by Goldstein and Schlag \cite{gs1} assuming that $\alpha$ is  strong Diophantine. In view of Avila's global theory \cite{avila0}, for one-frequency analytic quasi-periodic Schr\"odinger operators,  the energy can be divided into three  regimes: the subcritical regime (almost reducible regime), the supercritical regime (positive Lyapunov exponent regime) and the critical regime (otherwise), and typically, there is no critical energy \cite{avila0}. Thus the  H\"older continuity of IDS  for  one-frequency quasi-periodic Schr\"odinger operators is clear: it is  H\"older continuous in both subcritical regime and supercritical regime, while  one should not expect any modulus of continuity of IDS at the critical energies \cite{B1}. We remark that for small potentials,   \cite{amor,aj} got the optimal result, i.e., $\frac{1}{2}$-H\"older continuity, while for big potentials, \cite{gs1} did not give any information on the  H\"older exponent. So it is natural to ask: what is the  modulus of continuity of IDS for big potentials?   For the almost Mathieu operator, Bourgain \cite{b1} proved that for Diophantine $\alpha$ and large enough $\lambda$, the H\"older exponent is larger than $\frac{1}{2}-\e$ for any $\e>0$. The H\"older exponent  was  further optimized  as   $\frac{1}{2}$ by Avila and Jitomirskaya \cite{aj} for non-critical almost Mathieu operator with Diophantine frequency. For general 
trigonometric polynomial of degree $m$,  Goldstein and Schlag \cite{gs2} proved that the H\"older exponent is larger than $\frac{1}{2m}-\e$ for all $\e>0$ in the positive Lyapunov exponent regime, which is a generalization of  Bourgain's result \cite{b1}. 
Their methods are  the Large Deviation Theorem (LDT) and Avalanch Principle (AP), which were initially developed by Bourgain-Goldstein \cite{bg} and Goldstein-Schlag \cite{gs1}, further developed by Goldstein-Schlag \cite{gs2}. LDT and AP were later recognized as powerful tools to study the regularity of IDS for quasi-periodic Schr\"odinger operators in the positive Lyapunov exponent regime. Based on LDT and AP, various H\"older continuity results of IDS for quasi-periodic Schr\"odinger operators were obtained, here is a partial list \cite{b1,gs1,gs2,hz,LWY,wz1,YZ,xgw}. Initially, LDT and AP were developed for Schr\"odinger cocycles, which works only for $SL(2,\R)$ cocycles. More recently, LDT and AP were further developed to higher dimensional linear cocycles \cite{s,dk1,dk2} for proving the continuity of the Lyapunov exponents. We remark that the characterization of the H\"older exponent  is much more difficult in the supercritical regime, which is closely related to the  ``multiple resonances"  \footnote{Roughly speaking, on the almost reducibility side, resonance means the eigenvalues of the constant matrix come very close to each other. On the large deviation side, resonance means the eigenvalues of the operators restricted to a box come very close to each other.} and 
developing techniques to deal with  such ``multiple resonances" is very difficult (\cite{gs2} is an 114 pages paper). Different from the small potentials,  $\frac1{2}$-H\"older continuity of IDS usually  should not be expected for large potentials. It seems that  the H\"older exponent is closely related to the profile of the potential.  In \cite{gs2}, Goldstein and Schlag  linked 
the H\"older exponent to the degree of the trigonometric polynomial potential.

 In this paper, different from \cite{gs2}, we use the Quantitative Almost Reducibility Theorem (QART) and Aubry duality to study this problem. Our approach has some advantages: it does not sensitive to $d$ and $m$ (\cite{gs2} treated the case $m=d=1$), thus it works for more general classes of operators (\ref{1.2}) and (\ref{schrodinger-2}).   More importantly, it brings a more delicate estimate on the H\"older exponent. More precisely, for large trigonometric polynomial potentials, the H\"older exponent of IDS is controlled by the cardinality of the level sets $\#\{\theta: W(\theta)=E\}$ instead of  the degree of $W$.  Our main result is

\begin{Theorem}\label{main}
Assume that $\alpha\in {\rm DC}_d$, $\e>0$ and $\#\{\theta: W(\theta)=E\}$, i.e.,  the number of zeros (counting multiplicities) of $W(\theta)-E$ on $\T$, is no more than $2m_0$ for all $E\in \R$. Then there exist $\lambda_0=\lambda_{0}(\alpha,W,V,\e)>0$ and $C=C(\alpha,W,V)>0$ such that if $\lambda\geq \lambda_0$, for operators \eqref{1.2} and \eqref{schrodinger-2}, we have
$$
|N(E)-N(E')|\leq C|E-E'|^{\frac{1}{2m_0}-\e},
$$
for all $E,E'\in\R$.
\end{Theorem}
\begin{Remark}
We indeed characterize the H\"older exponent in the neighborhood of  $E$ by the cardinality of level sets $\#\{\theta: W(\theta)=E\}$. For $W$ in the following figure
\begin{figure}[htbp]
\includegraphics[width=0.4\linewidth]{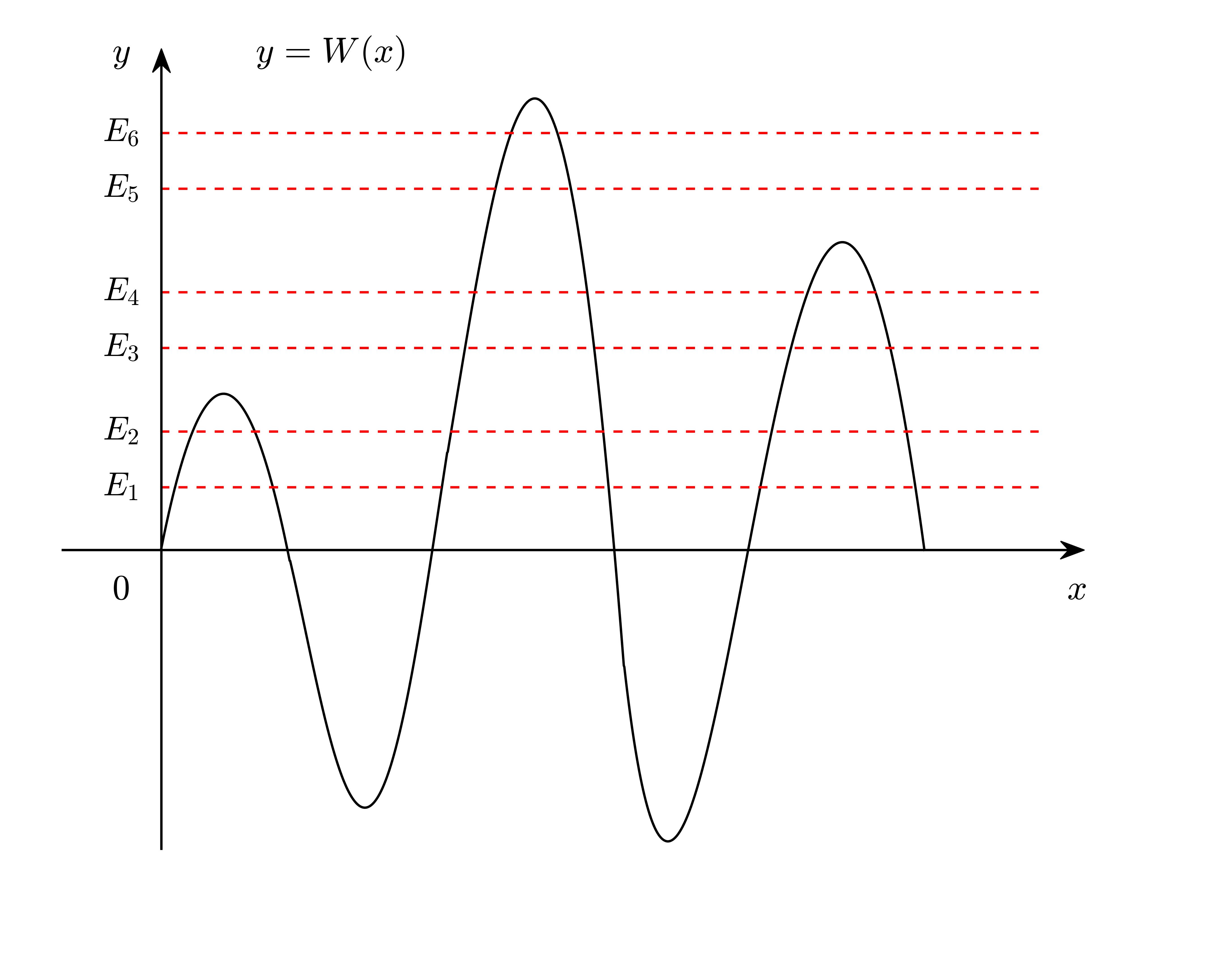}
\end{figure}
\\
we have 
\begin{enumerate}
\item For $E\in [E_1,E_2]$, the H\"older exponent is larger than $\frac{1}{6}-\e$,
\item For $E\in [E_3,E_4]$, the H\"older exponent is larger than  $\frac{1}{4}-\e$,
\item For $E\in [E_5,E_6]$,  the H\"older exponent is larger than  $\frac{1}{2}-\e$.
\end{enumerate}
\end{Remark}

\begin{Remark}
If $d=1$ and $V(x)=2\cos2\pi x$, it is possible to extend the result to operator (\ref{1.2}) with any real analytic potential $W$ and large $\lambda$. We  guess that the H\"older exponent of IDS is closely related to the acceleration \footnote{See section 2.4 for the definition.} in the positive Lyapunov exponent regime. 
\end{Remark}
\begin{Remark} In the perturbative regime, if $m=d=1$, our result  is a refinement of  Goldstein and Schlag's result \cite{gs2}. 
\end{Remark}
\begin{Remark}
To the best of our knowledge, Theorem \ref{main} is the first H\"older regularity result of IDS for quasi-periodic long-range operators (\ref{1.2}) on $\ell^2(\Z^d)$ with potentials beyond the cosine function.
It is also the first H\"older regularity result of IDS for multi-frequency finite difference analytic quasi-periodic operators (\ref{schrodinger-2}).
\end{Remark}
\begin{Remark}
 It is  possible to improve the H\"older exponent to $\frac{1}{2m_0}$ in Theorem \ref{main} by refining estimates of QART for $GL(m,\C)$ cocycles. We do not go further in the present paper.
\end{Remark}

Theorem \ref{main} is a generalization of \cite{gs2} from one dimensional Schr\"odinger operators to Schr\"odinger operators on $\Z^d$ lattice in the perturbative regime. In fact, since $m_0\leq m$, we immediately get the following  corollary which was proved by Goldstein and Schlag in \cite{gs2} for the case $d=1$ and $V(x)=2\cos2\pi x$.
\begin{Corollary}
Assume $\alpha\in {\rm DC}_d$ and $\e>0$, there exist $\lambda_0=\lambda_{0}(\alpha,W,V,\e)$ and $C_0=C(\alpha,W,V)$ such that if $\lambda\geq \lambda_0$, for operators \eqref{1.2} and \eqref{schrodinger-2}, we have
$$
|N(E)-N(E')|\leq C|E-E'|^{\frac{1}{2m}-\e},
$$
for all $E,E'\in\R$.
\end{Corollary}

We briefly review some other related regularity results on IDS. The H\"older continuity of IDS for one-frequency analytic quasi-periodic Schr\"odinger operators with Liouvillean frequencies was obtained by You-Zhang \cite{YZ} and  Han-Zhang \cite{hz}.  Weak H\"older continuity of IDS for quasi-periodic Schr\"odinger operators on $\ell^2(\Z^2)$ with Diophantine frequencies and large analytic potentials was proved by Schlag \cite{schlag}. For the lower regularity case, Klein \cite{klein} proved that for Schr\"{o}dinger operators with potentials in a Gevrey class, the IDS  is  weak H\"{o}lder continuous on any compact interval of the energy provided that the coupling constant is large enough, the frequency is Diophantine and the potential satisfies some transversality condition. Wang-Zhang \cite{wz1} obtained the weak H\"{o}lder continuity of the IDS as function of energies, for a class of $C^2$ quasi-periodic potentials and for any Diophantine frequency. Subsequently it was improved to be H\"older continuous  by Liang-Wang-You \cite{LWY}. Recently, the H\"older exponent was proved to be $1/2$ by Xu-Ge-Wang \cite{xgw} which is optimal. More recently, Cai-Chavaudret-You-Zhou \cite{ccyz} proved $1/2$-H\"older continuity of IDS for quasi-periodic Schr\"odinger operator with small finitely differential potentials and Diophantine frequencies. Jitomirskaya-Kachkovskiy \cite{jk} proved Lipschitz continuity of IDS for quasi-periodic Schr\"odinger operators with bounded monotonic potentials. Their result \cite{jk} was recently extended by Kachkovskiy \cite{kach} to unbounded monotonic potentials.
\subsection{Regularity of the Lyapunov exponents  for quasi-periodic $GL(m,\C)$ cocycles} 
We denote by $GL(m,\C)$ the set of all $m\times m$ invertible matrices. Given $A \in C^\omega(\T^d, GL(m,\C))$ and rational independent $\alpha\in\R^d$, we define the {\it quasi-periodic $GL(m,\C)$ cocycle} $(\alpha,A)$:
$$
(\alpha,A)\colon \left\{
\begin{array}{rcl}
\T^d \times \C^{m} &\to& \T^d \times \C^{m}\\[1mm]
(x,v) &\mapsto& (x+\alpha,A(x)\cdot v)
\end{array}
\right.  .
$$
The iterates of $(\alpha,A)$ are of the form $(\alpha,A)^n=(n\alpha, A_n)$, where
$$
A_n(x):=
\left\{\begin{array}{l l}
A(x+(n-1)\alpha) \cdots A(x+\alpha) A(x),  & n\geq 0\\[1mm]
A^{-1}(x+n\alpha) A^{-1}(x+(n+1)\alpha) \cdots A^{-1}(x-\alpha), & n <0
\end{array}\right.    .
$$
We denote by $L_1(\alpha, A)\geq L_2(\alpha,A)\geq...\geq L_m(\alpha,A)$ the Lyapunov exponents of $(\alpha,A)$ repeatedly according to their multiplicities, i.e.,
$$
L_k(\alpha,A)=\lim\limits_{n\rightarrow\infty}\frac{1}{n}\int_{\T^d}\ln\sigma_k(A_n(x))dx,
$$
where for a matrix $B\in GL(m,\C)$, we denote by $\sigma_1(B)\geq...\geq \sigma_m(B)$ its singular values (eigenvalues of $\sqrt{B^*B}$). 

The continuity of the Lyapunov exponents of linear cocycles has been  extensively studied. It was proved by Bourgain-Jitomirskaya  in \cite{bj} that LE is joint continuous for $SL(2,\R)$ cocycles, in frequency and cocycle map, at any irrational frequencies.  Jitomirskaya-Koslover-Schulteis \cite{jks} got the continuity of LE with respect to potentials for a class of
analytic quasi-periodic $M(2,\C)$ cocycles.  Bourgain \cite{bourgain2} extended the results in \cite{bj} to multi-frequency case. 
Jitomirskaya-Marx \cite{jm} extended the results in \cite{bj} to all (including singular) $M(2,\C)$ cocycles. More recently, continuity of the Lyapunov exponents for one-frequency analytic $M(m,\C)$ cocycles was given by Avila-Jitomirskaya-Sadel \cite{ajs}. Weak H\"older continuity of the Lyapunov exponents for multi-frequency $GL(m, \C)$ cocycles, $m \geq 2,$ was recently obtained by Schlag \cite{s} and
Duarte-Klein \cite{dk1}.

In this paper, we show that if $A\in C^\omega(\T^{d},GL(m,\C))$, the Lyapunov exponents can still be H\"older continuous, provided that the cocycle is almost reducible. Our main result is the following:
            
\begin{Theorem}\label{thm1}
Let $\alpha \in DC_d(\kappa,\tau)$ and $\e>0$, if $(\alpha,A)$ is $C^\omega_h$ almost reducible, then for any $\widetilde{A}\in C^\omega(\T^d,GL(m,\C))$, we have
$$\lvert L_i(\alpha,A)-L_i(\alpha,\widetilde{A})\rvert \leq C_0 |\widetilde{A}-A|_0^{\frac{1}{2n_i}-\e},\ \ 1\leq i\leq m,$$
where $C_0$ is a constant depending on $A,d,\kappa,\tau,m,\e$ and $n_i$ is the multiplicity of $L_i(\alpha,A)$.
\end{Theorem}
\begin{Remark}
Almost reducible cocycles include a large class of cocycles. For example, consider $(\alpha, A_0e^{f_0(\cdot)})$ where $A_0\in GL(m,\C)$ and $f_0\in C_h^\omega(\T^d,gl(m,\C))$, then $(\alpha, A_0e^{f_0(\cdot)})$ is almost reducible provided $|f_0|_h$ is sufficiently small. In one-frequency case, all subcritical cocycles are almost reducible \cite{avila1,avila2}.
\end{Remark}

\subsection {Outline of the proofs.} 
The Quantitative Almost Reducibility Theorem (QART), was initially developed by Eliasson \cite{e}, further developed by Hou-You \cite{hy}, Avila-Jitomirskaya \cite{aj} and Avila \cite{avila1,avila2}. It has been proved to be a powerful tool  in studying various spectral problems of quasi-periodic operators in the almost reducible regime which is the ``dual regime" of the positive Lyapunov exponent regime \cite{avila0} \footnote{For example, for the almost Mathieu operator, the positive Lyapunov exponent regime corresponds to $|\lambda|>1$ and the almost reducible regime corresponds to $|\lambda|<1$.}.

We will prove our main results  by establishing  the QART for higher dimensional linear cocycles. More precisely, we establish a quantitative version of the almost reducibility theorem in the higher dimension emphasizing the quantitative estimates on the Lyapunov exponents at each KAM iteration. The advantage of QART is that it works for higher dimensional cocycles with multiple frequencies, which seems to be highly nontrivial  by LDT and AP.  Let us explain why QART is a powerful tool for investigating the regularity of the Lyapunov exponents of analytic quasi-periodic $GL(m,\C)$ cocycles.
\begin{definition}
$(\alpha,A)\in C^\omega(\T^d,GL(m,\C))$ is  said to be $C^\omega_h$ almost reducible  if there exist $B_j\in C_h^\omega(\T^d,GL(m,\C))$, $A_j\in GL(m,\C)$ and $f_j\in C_h^\omega(\T^d,gl(m,\C))$ such that
$$
B_j(x+\alpha)A(x)B_j^{-1}(x)=A_je^{f_j(x)},
$$
with $|f_j|_h\rightarrow 0$ and $A_j\rightarrow A_{\infty} \in GL(m,\C)$. 
\end{definition}

\noindent The followings are three basic facts of the Lyapunov exponents:
\begin{enumerate}
\item The Lyapunov exponents are invariant under continuous conjugation, i.e.,
$$
L_i(\alpha,A)=L_i(\alpha,B_j(\cdot+\alpha)A(\cdot)B_j^{-1}(\cdot))=L_i(\alpha,A_je^{f_j(\cdot)}),\ 1\leq i\leq m.
$$
Furthermore, if the base frequency is Diophantine, all Lyapunov exponents are continuous in $A$ in $C^\omega_h$ topology \cite{s,dk2}.  It follows that
$$
L_i(\alpha,A_je^{f_j(\cdot)})\rightarrow  L_i(\alpha,A_{\infty}),\ 1\leq i\leq m.
$$
Thus,
$$
L_i(\alpha,A)= L_i(\alpha,A_{\infty}), \  1\leq i\leq m.
$$
\item Almost reducibility is open \footnote{It follows from the results in \cite{eliasson,chavaudret}.}. Thus if $|\widetilde{A}-A|_h$ is sufficiently small, $(\alpha,\widetilde{A})$ is also almost reducible to $(\alpha,\widetilde{A}_\infty)$ with
$$
L_i(\alpha,\widetilde{A})=L_i(\alpha,\widetilde{A}_{\infty}),\ 1\leq i\leq m.
$$
\item The Lyapunov exponents of constant cocycles are computable. If we denote the eigenvalues of $A_{\infty}$ ($\widetilde{A}_{\infty}$) by $\{e^{-2\pi i \rho_i}\}_{i=1}^m$ ($\{e^{-2\pi i {\tilde{\rho}}_i}\}_{i=1}^m$) respectively, then
$$
\{L_i(\alpha,A)\}_{i=1}^m=\{2\pi\Im \rho_i\}_{i=1}^m,\ \  \{L_i(\alpha,\widetilde{A})\}_{i=1}^m=\{2\pi\Im \tilde{\rho}_i\}_{i=1}^m.
$$
\end{enumerate}
QART can give us very precise estimates on the differences $|L_i(\alpha,A_\infty)-L_i(\alpha,\widetilde{A}_\infty)|$, from which one can find H\"older continuity easily. The H\"older continuity of IDS is then a consequence in view of the generalized Thouless formula \cite{craigsimon1,craigsimon2}.

Finally, we introduce our motivations. During the past ten years, QART for $SL(2,\R)$ cocycles along with Aubry duality \cite{gjls} has been proved to be a powerful tool to solve various central problems in the field of spectral theory of quasi-periodic operators, see \cite{ayz1,ayz2,aj,jk1,gyz1,gyz,gy} for some recent progresses. While all the results are essentially restricted to the quasi-periodic long-range operator on $\ell^2(\Z^d)$ with a cosine potential or the almost Mathieu operator. When the potential goes beyond cosine, things become dramatically complicated, and the so-called ``multiple resonances" becomes a real issue.  KAM theory provides us a way to deal with ``multiple resonances" \cite{YJDE} and it may have wider applications in other aspects of the spectral theory for quasi-periodic Schr\"odinger operators with potentials beyond the cosine function.

\section{Preliminaries}
For a bounded
analytic (possibly matrix valued) function $F$ defined on $ \{x |  | \Im x |< h \}$, let
$
|F|_h=  \sup_{ | \Im x |< h } | F(x)| $ and denote by $C^\omega_{h}(\T^d,*)$ the
set of all these $ *$-valued functions ($*$ will usually denote $\R$, $\C$, $gl(m,\C)$ and
$GL(m,\C)$). Also we denote $C^\omega(\T^d,*)=\cup_{h>0}C^\omega_{h}(\T^d,*)$, and $|F| _{0}=\sup_{x \in \T^{d}}|F(x)|.$
\subsection{Quasi-periodic finite-range cocycles}
In this subsection, we review some basic concepts for quasi-periodic finite-range operators, especially for quasi-periodic Schr\"odinger operators. An important example of quasi-periodic $GL(2m,\C)$ cocycle is the finite-range cocycle  $(\alpha,L_{E,W}^{\lambda^{-1}V})$ where
$$
L_{E,W}^{\lambda^{-1}V}(x)=\frac{1}{\widehat{W}_m}
\begin{pmatrix}
\begin{smallmatrix}
-\widehat{W}_{m-1}&\cdots&-\widehat{W}_1&E-\lambda^{-1}V(x)-\widehat{W}_0&-\widehat{W}_{-1}&\cdots&-\widehat{W}_{-m+1}&-\widehat{W}_{-m}\\
\widehat{W}_m& \\
& &  \\
& & & \\
\\
\\
& & &\ddots&\\
\\
\\
& & & & \\
& & & & & \\
& & & & & &\widehat{W}_{m}&
\end{smallmatrix}
\end{pmatrix}.
$$
The finite-range cocycles are equivalent to the eigenvalue equations of operators \eqref{schrodinger-2}, i.e., $L_{W,\alpha,x}^{\lambda^{-1}V}u=E u$. Note that $(\alpha,L_{E,W}^{\lambda^{-1}V})$ is more special than general $GL(2m,\C)$ cocycles in the following senses,
\begin{enumerate}
\item The $m$-th iteration of $(\alpha,L_{E,W}^{\lambda^{-1}V})$, denoted by $(m\alpha,(L_{E,W}^{\lambda^{-1}V})_m)$,  is a symplectic cocycle \cite{haro}. As a corollary,  the Lyapunov exponents of  $(\alpha,L_{E,W}^{\lambda^{-1}V})$ come into pairs $\pm L_i(\alpha,L_{E,W}^{\lambda^{-1}V})$ ($1\leq i\leq m$).
\item We denote $L_{i}(E)=L_i(\alpha,L_{E,W}^{\lambda^{-1}V})\geq 0$ ($1\leq i\leq m$) and the IDS of $(L^{\lambda^{-1} V} _{W,\alpha,x})_{x\in\T^d}$ by $N(E)$, then the sum of the nonnegative Lyapunov exponents and the IDS are linked by the famous generalized Thouless formula \cite{craigsimon1,craigsimon2,haro},
\begin{equation}
\sum\limits_{i=1}^m L_i(E)+\ln|\widehat{W}_m|=\int\ln|E-E'|dN(E').
\end{equation}
\end{enumerate}
By the Hilbert transform and the theory of singular integral operators, the H\"older continuity passes from $\sum\limits_{i=1}^m L_i(E)$ to $N(E)$ and vice versa (see \cite{gs1} for details). Notice that the Aubry duals of $(L^{\lambda^{-1} V} _{W,\alpha,x})_{x\in\T^d}$ are $(L^{W} _{\lambda^{-1}V,\alpha,\theta})_{\theta\in\T}$, we denote the IDS of  $(L^{W} _{\lambda^{-1}V,\alpha,\theta})_{\theta\in\T}$ by $\widehat{N}(E)$. The IDS is invariant under Aubry dual, i.e.,
\begin{Proposition}[\cite{haro,jm1}]
$N(E)=\widehat{N}(E)$.
\end{Proposition}
\subsection{Some basic properties of the Lyapunov exponents}
In the introduction, we have given some basic facts on the Lyapunov exponents of linear cocycles. In this subsection, we give an elementary proof of them. \\
1. Lyapunov exponents are invariant under continuous conjugations.
\begin{Proposition}\label{invariance}
Assume $(\alpha,A)\in \T^d\times C^0(\T^d,GL(m,\C)),\ \ B\in C^0(\T^d,GL(m,\C)),$ and $\widetilde{A}(x)=B(x+\alpha)A(x)B^{-1}(x)$, we have
$$
L_i(\alpha,\widetilde{A})=L_i(\alpha,A),\ \ 1\leq i\leq m.
$$
\end{Proposition}
\begin{pf}
For any $1\leq i\leq m$, we have
$$
\Lambda^i\widetilde{A}_n(x)=\Lambda^iB(x+n\alpha)A_n(x)B^{-1}(x),
$$
where $\Lambda^i M$ is the $i$-th wedge of $M$. Thus
$$
|B|_0^{-m}|B^{-1}|_0^{-m}|\Lambda^iA_n(x)|\leq |\Lambda^i\widetilde{A}_n(x)|\leq |B|_0^{m}|B^{-1}|_0^m|\Lambda^iA_n(x)|,
$$
which implies that
$$
\frac{1}{n}\int_{\T^d}\ln|\Lambda^i\widetilde{A}_n(x)|dx=\frac{1}{n}\int_{\T^d}\ln|\Lambda^iA_n(x)|dx+o(1).
$$
Hence
$$
L_i(\alpha,\widetilde{A})=L_i(\alpha,A),\ \ 1\leq i\leq m.
$$
\end{pf}

\noindent 2. The Lyapunov exponents of constant cocycles are computable. 
\begin{Proposition}\label{leconstant}
If we denote the eigenvalues of $A\in GL(m,\C)$ by $\{e^{-2\pi i \rho_j}\}_{j=1}^m$, then
$$
\{L_j(\alpha,A)\}_{j=1}^m=\{2\pi\Im \rho_j\}_{j=1}^m.
$$
\end{Proposition}
\begin{pf}
We only give the proof of $L_1(\alpha,A)=\max_{1\leq j\leq m}\{2\pi\Im \rho_j\}$ and the other proofs are similar. By Proposition \ref{invariance}, we can always assume that $A$ is in Jordan form. Then one can directly compute that
$$
L_1(\alpha,A)=\lim\limits_{n\rightarrow\infty}\frac{1}{n}\int_{\T^d}\ln\sigma_1(A^n)dx=\lim\limits_{n\rightarrow\infty}\frac{1}{n}\ln\sigma_1(A^n)=\max_{1\leq j\leq m}\{2\pi\Im \rho_j\}.
$$

\end{pf}
\subsection{Perturbation theory of constant matrices}
In this subsection, we briefly introduce the perturbation theory of constant matrices
which will be used in Section 3.
\begin{Definition} Let $A$ be an $m\times m$ matrix, denote
$$
|A|_F=\sqrt{\sum\limits_{i=1}^m\sum\limits_{j=1}^m|a_{ij}|^2}.
$$
\end{Definition}
We call $A$ normal if $AA^T=A^TA$ where $A^T$ is the transpose conjugation of
$A$.
\begin{Lemma}[\cite{song}]\label{per}
Let $A$ and $B$ be two $m\times m$ matrices, where $A$ is normal and
$B$ is nonnormal, with spectrum
$\lambda(A)=\{\lambda_1,\cdots,\lambda_m\}$ and
$\lambda(B)=\{\mu_1,\cdots,\mu_m\}$, then there exists a
permutation $\pi$ of $\{1,2,\cdots,m\}$ such that
$$
\sqrt{\sum\limits_{j=1}^m\left|\mu_{\pi(j)}-\lambda_j\right|^2}\leq\sqrt{m}|B-A|_F.
$$
\end{Lemma}
The following proposition follows immediately.
\begin{Proposition}\label{eigen_per}
Let $A$ and $B$ be two $m\times m$ matrices, with spectrum
$\lambda(A)=\{\lambda_1,\cdots,\lambda_m\}$ and
$\lambda(B)=\{\mu_1,\cdots,\mu_m\}$, then there exists a
permutation $\pi$ of $\{1,2,\cdots,m\}$ such that
$$
\sqrt{\sum\limits_{j=1}^m\left|\mu_{\pi(j)}-\lambda_j\right|^2}\leq C(m,A)|B-A|^{\frac{1}{m}}.
$$
\end{Proposition}
\begin{pf}
Let $\epsilon=|B-A|$. There exists a unitary matrix $U\in GL(m,\mathbb{C})$ such that
$$
U^{-1}AU=\begin{pmatrix}\lambda_1&*&*\\&\ddots&* \\&&\lambda_m\end{pmatrix}.
$$
If we denote $T=diag\{1,\epsilon^{\frac{1}{m}},\cdots, \epsilon^{\frac{m-1}{m}}\}$, then
\begin{align*}
(UT)^{-1}BUT&=(UT)^{-1}AUT+(UT)^{-1}(B-A)UT\\
&=diag\{\lambda_1,\cdots,\lambda_m\}+F,
\end{align*}
with $|F|\leq 2|A|\epsilon^{\frac{1}{m}}$. By Lemma \ref{per}, there exists a
permutation $\pi$ of $\{1,2,\cdots,m\}$ such that
$$
\sqrt{\sum\limits_{j=1}^m(\mu_{\pi(j)}-\lambda_j)^2}\leq \sqrt{m}|F|_F\leq C(m,A)|B-A|^{\frac{1}{m}}.
$$
\end{pf}

\subsection{Global theory for analytic one-frequency quasi-periodic Schr\"odinger operators}
In 2015,  Avila \cite{avila0} gave a qualitative spectral
picture for one-frequency quasi-periodic Schr\"odinger operators. To explain more, we denote the associated cocycle of the eigenequation $H_{W(\cdot+i\e),\alpha,\theta}u=Eu$ by $(\alpha, S_{E}^W(\cdot+i\e))$ and the associated nonnegative Lyapunov exponent  by $L_\e(E)$. The {\it acceleration} is defined in \cite{avila0} by
$$
\omega(E)=\lim\limits_{\e\rightarrow 0^+}\frac{L_\e(E)-L(E)}{2\pi\e}.
$$
The global theory discovers that the spectral set of the Schr\"odinger operator can be divided into three regimes based on the Lyapunov exponent and acceleration:
\begin{enumerate}
\item The subcritical regime: $L(E)=0$ and $\omega(E)=0$.
\item The critical regime: $L(E)=0$ and $\omega(E)>0$.
\item The supercritical regime: $L(E)>0$ and $\omega(E)>0$.
\end{enumerate}
Moreover, the subcritical regime is equivalent to the almost reducible regime, announced by Avila \cite{avila1,avila2}. Usually, the almost reducible regime can be viewed as the dual regime of the positive Lyapunov exponent regime.

\section{A quantitative almost reducibility theorem for $GL(m,\C)$ cocycles }

Consider the local quasi-periodic $GL(m,\C)$ cocycle
\begin{align*}
(\alpha,A_0e^{f_0(\cdot)})\colon \left\{
\begin{array}{rcl}
\T^d \times \C^{m} &\to& \T^d \times \C^{m}\\[1mm]
(x,v) &\mapsto& (x+\alpha,A_0e^{f_0(x)}\cdot v)
\end{array}
\right.  ,
\end{align*}
where $A_0\in GL(m,\C)$ is a constant matrix,  $f_0\in C_h^\omega(\T^d, gl(m,\C))$ is an analytic perturbation.  $\alpha\in DC_d(\kappa,\tau)$ for some $\kappa>0$, $\tau>d-1$.
In this section,  we prove the following  quantitative almost reducibility result. For our purpose, we specially pay attention to shift of the norms of the eigenvalues during the KAM iteration. 

\begin{Theorem}\label{thm3}
For any given $0<\sigma<\frac{1}{500m^3}$ and $0<h'<h$, there exists $\epsilon_0=\epsilon(\alpha,h,h',m, |A_0|,\sigma)$ such that if $|f_0|_h\leq \epsilon_0$, then $\left(\alpha , A_0e^{f_0(\cdot)}\right)$ is $C_{h'}^\omega$ almost reducible, i.e., there exist $B_{j}\in C^{\omega}_{h'}(\T^{d}$,
$GL(m,\C))$, $A_{j}\in GL(m,\C)$ and $f_{j}\in C^{\omega}_{h'}(\T^d,gl(m,\C))$, such that
$$
B_{j}(x+\alpha)(A_0e^{f_0(x)})B^{-1}_{j}(x)=A_{j}e^{f_j(x)}.
$$
Moreover, we have the following estimates
\begin{equation}\label{ess3}
|f_{j}|_{h'}\leq \epsilon_{j},\ \ |B_{j}|_0\leq \epsilon_{j-1}^{-200m^2\sigma},
\end{equation}
\begin{align}\label{ess2}
\left|\sum\limits_{\ell=1}^m\Im\rho_{\ell}^j-\sum\limits_{\ell=1}^m\Im\rho_{\ell}^{j-1}\right|\leq C(m,A_0)\epsilon_{j-1}^{\frac{1}{2}},
\end{align}
\begin{equation}\label{ess4}
\left| \Im\rho_{\ell}^j-\Im\rho_{\ell}^{j-1}\right|\leq C(m,A_0)\epsilon^{\frac{1}{m}-200m^2\sigma}_{j-1},\ \ 1\leq \ell\leq m,
\end{equation}
where $\left\{e^{-2\pi i\rho_\ell^j}\right\}_{\ell=1}^m$ are the eigenvalues of $A_j$ and $\epsilon_j= \epsilon_0^{2^j}$.
\end{Theorem}
\pf Suppose that
\begin{align}\label{initial_assump}
|f_0|_h\leq \epsilon_0 \leq \frac{C}{|A_0|^{\frac{500dm^3\tau}{\sigma}}}(h-h')^{\frac{500dm^3\tau}{\sigma}},
\end{align}
where $C=\min\left\{2^{-\frac{2\tau}{\sigma}}(4m+1)^{-\frac{2}{\sigma}}\kappa^{\frac{2}{\sigma}},C_0^{-\frac{500dm^3\tau}{\sigma}}\right\}$ \footnote{$C_0$ is an absolute constant only depending on $m$.}. Then we define the following sequences inductively,
$$
\epsilon_{j+1}=\epsilon_{j}^{2},\ \ h_j-h_{j+1}=\frac{h-h'}{4^{j+1}}, \ \ N_j=\frac{2|\ln\epsilon_j|}{h_j-h_{j+1}}.
$$
By our selection of  $\epsilon_0$, it's easy to check that
\begin{equation}\label{iter}
\epsilon_j \leq \frac{C}{|A_j|^{\frac{500dm^3\tau}{\sigma}}}(h_j-h_{j+1})^{\frac{500dm^3\tau}{\sigma}}.
\end{equation}
Indeed, $\epsilon_j$ on the left side of the inequality decays at least super-exponentially with $j$, while $(h_j-h_{j+1})^{\frac{500dm^3\tau}{\sigma}}$ on the right side decays exponentially with $j$.

Assume after $j$ steps of iteration, we are at the $(j+1)^{th}$ KAM step. That is, we have already constructed $B_j\in C^\omega_{h_{j}}(\T^d,GL(m,\C))$ such that
\begin{equation}\label{eq_1}
B_{j}(x+\alpha)A_0e^{f_0(x)}B^{-1}_{j}(x)=A_{j}e^{f_{j}(x)},
\end{equation}
where $A_j\in GL(m,\C)$ with eigenvalues $\{e^{-2\pi i\rho^{j}_\ell}\}_{\ell=1}^m$ and
\begin{equation}\label{ess3'}
|f_{j}|_{h'}\leq \epsilon_{j},\ \ |B_{j}|_0\leq \epsilon_{j-1}^{-200m^2\sigma},
\end{equation}
\begin{align*}
\left|\sum\limits_{\ell=1}^m\Im\rho_{\ell}^j-\sum\limits_{\ell=1}^m\Im\rho_{\ell}^{j-1}\right|\leq C(m,A_0)\epsilon_{j-1}^{\frac{1}{2}},
\end{align*}
\begin{equation*}
\left| \Im\rho_{\ell}^j-\Im\rho_{\ell}^{j-1}\right|\leq C(m,A_0)\epsilon^{\frac{1}{m}-200m^2\sigma}_{j-1},\ \ 1\leq \ell\leq m.
\end{equation*}
We will construct
$$
\bar{B}_j\in C^\omega_{h_{j+1}}(\T^d,GL(m,\C)),\ \ A_{j+1}\in GL(m,\C),\ \ f_{j+1}\in C_{h_{j+1}}(\T^d,gl(m,\C))
$$
such that
$$
\bar{B}_j(x+\alpha)A_je^{f_j(x)}\bar{B}^{-1}_j(x)=A_{j+1}e^{f_{j+1}(x)},
$$
with desired estimates. The proof is divided into the following four steps. We denote $A_j, f_j,\epsilon_j,h_j,\epsilon_{j+1},h_{j+1}$ by $A, f,\epsilon,h,\epsilon_+,h_+$ for simplicity.\\

\noindent
{\bf Step 1: Block diagonalizing $A$.} We denote $\{e^{-2\pi i\rho_\ell}\}_{1\leq \ell\leq m}$ the eigenvalues of $A$.
\begin{Lemma}[Block diagonalization]\label{class}
$\{\rho_\ell\}_{1\leq \ell\leq m}$ can be grouped as  $\bigcup\limits_{\ell=1}^r E_{\ell}$ with $\# E_{\ell}=n_\ell$, such that
\begin{equation}\label{g1}
\left|\rho-\rho'\right|\leq m\epsilon^{\sigma},\ \ \rho,\rho'\in E_{\ell},
\end{equation}
\begin{equation}\label{g2}
\left|\rho-\rho'\right|\geq \epsilon^{\sigma},\ \ \rho\in E_{k},\rho'\in E_{\ell},\ \ k\neq \ell.
\end{equation}
Moreover, there exist $C_0\geq 1$ depending only on $m$, $P\in GL(m,\C)$ with estimate
\begin{equation}\label{g3}
|P^{-1}|, \ \ |P| \leq C_0(2|A|\epsilon^{-\sigma})^{m(m+1)},
\end{equation}
and  upper triangular matrices $\{A_{\ell}\}_{\ell=1}^r$ with $spec(A_{\ell})=\{e^{-2\pi i \rho}|\rho\in E_{\ell}\}$, such that
$$
PAP^{-1}=diag\{A_{1}\cdots A_{r}\}.
$$

\end{Lemma}
\begin{pf}
\eqref{g1} and \eqref{g2} follow from a simple observation. \eqref{g3} follows from Lemma $A^1$ in \cite{eliasson}.
\end{pf}
Thus we can conjugate $\left(\alpha, Ae^{f(\cdot)}\right)$ to a new cocycle,
\begin{equation}\label{new1}
PAe^{f(x)}P^{-1}=diag\{A_{1}\cdots A_{r}\}e^{Pf(x)P^{-1}}:=\widetilde{A}e^{\tilde{f}(x)}.
\end{equation}
Since $C\leq C_0^{-\frac{500dm^3\tau}{\sigma}}$, together with the assumption that $\sigma<\frac{1}{500m^3}$, we have
\begin{equation}\label{new_f}
|\tilde{f}|_h \leq |P||f|_h|P^{-1}|\leq \epsilon C^2_0(2|A|\epsilon^{-\sigma})^{2m(m+1)}\leq \min \left\{\epsilon^{1-200m^2\sigma},\epsilon^{\frac{9}{10}}\right\}.
\end{equation}
\\
\noindent
{\bf Step 2: Eliminating the non-resonant terms.}
For any given $\alpha\in \R^{d}$ and $A\in GL(m,\C)$, we decompose $\mathcal{B}_h=C^{\omega}_{h}(\T^{d},gl(m,\C))=\mathcal{B}_h^{nre}(\epsilon^{\frac{2}{5}}) \bigoplus\mathcal{B}_h^{re}(\epsilon^{\frac{2}{5}})$ in such a way that for any $Y\in\mathcal{B}_h^{nre}(\epsilon^{\frac{2}{5}})$,
\begin{equation}\label{space}
A^{-1}Y(\cdot+\alpha)A\in\mathcal{B}_h^{nre}(\epsilon^{\frac{2}{5}}), \qquad \lvert A^{-1}Y(\cdot+\alpha)A-Y(\cdot)\rvert_h\geq\epsilon^{\frac{2}{5}}\lvert Y\rvert_h.
\end{equation}
Now we define
\begin{align*}
\Lambda_{k,\ell}(\epsilon^{\sigma})=\left\{n\in\Z^{d}: \| E_{k}-E_{\ell}-\langle n,\alpha\rangle\|_{\R/\Z} \geq \epsilon^{\sigma}\right\}, \ \ 1\leq k,\ell\leq r. 
\end{align*}
where $\| E_{k}-E_{\ell}-\langle n,\alpha\rangle\|_{\R/\Z}=\min\limits_{\rho\in E_{k},\rho'\in E_{\ell}}\|\rho-\rho'-\langle n,\alpha\rangle\|_{\R/\Z}$.
\begin{align*}
\Lambda_N=\left\{f\in C^{\omega}_{h}(\T^{d},gl(m,\C))\mid f(x)=\sum\limits_{1\leq k,\ell\leq r}\sum_{n\in \Lambda_{k,\ell}(\epsilon^{\sigma})}\hat{f}_{k,\ell}(n)e^{2\pi i\langle n,x\rangle}\right\},
\end{align*}
where $\hat{f}_{k,\ell}(n)=P_{[\sum_{j=1}^{k-1}n_{j}+1,\sum_{j=1}^{k}n_j]}\hat{f}(n)P_{[\sum_{j=1}^{\ell-1}n_{j}+1,\sum_{j=1}^\ell n_j]}$ and $P_I:[1,m]\rightarrow I$ is a projection. The following lemma gives a characterization of the non-resonant space.
\begin{Lemma}\label{non-eli}
For any $Y\in \Lambda_N$, we have
$$
\left\lvert \widetilde{A}^{-1}Y(\cdot+\alpha)\widetilde{A}-Y(\cdot)\right\rvert_h\geq\epsilon^{\frac{2}{5}}\lvert Y\rvert_h.
$$
\end{Lemma}
\begin{pf}
Notice that $\widetilde{A}=diag\{A_{1},\cdots,A_{r}\}$. Thus for any $Y\in \Lambda_N$, we have
\begin{align*}
Y(x+\alpha)\widetilde{A}-\widetilde{A}Y(x)=\left(Y_{k,\ell}(x+\alpha)A_{\ell}-A_{k}Y_{k,\ell}(x)\right)_{1\leq k,\ell\leq r},
\end{align*}
where $Y_{k,\ell}(x)=P_{[\sum_{j=1}^{k-1}n_{j}+1,\sum_{j=1}^{k}n_j]}Y(x)P_{[\sum_{j=1}^{\ell-1}n_{j}+1,\sum_{j=1}^\ell n_j]}$. For any $(k,\ell)$ block, we apply Lemma \ref{useful_lem} by taking $A=A_{\ell}$, $B=A_{k}$, $\Lambda=\Lambda_{k,\ell}(\epsilon^{\sigma})$ and $\eta=\epsilon^{\sigma}$. By the definition of $\Lambda_N$ and \eqref{initial_assump}, it's easy to verify that condition \eqref{noncond} in Lemma \ref{useful_lem} is satisfied. Thus
\begin{align*}
\left\lvert \widetilde{A}^{-1}Y(\cdot+\alpha)\widetilde{A}-Y(\cdot)\right\rvert_h\geq \sum\limits_{1\leq k,\ell\leq r}\epsilon^{\sigma}\left(1+m\left(|A_{k}|+|A_{\ell}|\right)\epsilon^{-\sigma}\right)^{-2m}\left\lvert Y_{k,\ell}\right\rvert_h\geq\epsilon^{\frac{2}{5}}\lvert Y\rvert_h,
\end{align*}
the last inequality holds because of $\sigma<\frac{1}{500m^3}$ and \eqref{initial_assump}.
\end{pf}
This implies that $\Lambda_N\subset\mathcal{B}_h^{nre}(\epsilon^{\frac{2}{5}})$. Since $\epsilon^{\frac{2}{5}}\geq 20 |A|^2 \epsilon^{\frac{9}{20}}$,  we can apply  Lemma \ref{lem2} to remove all the non-resonant terms of $\tilde{f}$, which means there exist $Y\in \mathcal{B}_h$ and $\tilde{f}^{re}\in \mathcal{B}_h^{re}(\epsilon^{\frac{2}{5}})$ such that
\begin{equation}\label{new2}
e^{Y(x+\alpha)}\widetilde{A}e^{\tilde{f}(x)}e^{-Y(x)}=\widetilde{A}e^{\tilde{f}^{re}(x)},
\end{equation}
with $\lvert Y \rvert_h\leq \epsilon^{\frac{1}{3}}$ and $\lvert \tilde{f}^{re}\rvert_h\leq  2\min\left\{\epsilon^{1-200m^2\sigma},\epsilon^{\frac{9}{10}}\right\}$.\\

\noindent
{\bf Step 3:  Structure of the resonant terms  $\tilde{f}^{re}$.} By Lemma \ref{class}, we assume that 
$$
\tilde A=diag\{A_{1}\cdots A_{r}\}
$$ 
is a block diagonal matrix with $spec (A_{\ell})=\{e^{-2\pi i\rho}| \rho\in E_{\ell}\}$ ($\ell=1, \cdots, r$), satisfying (\ref{g1}) and (\ref{g2}). We say $A_{k}$ and  $A_{\ell}$ are $N$-resonant  if there are $\rho\in E_{k}$, $\rho'\in E_{\ell}$, $n_{k,\ell}\in \Z^d$ with $0<\lvert n_{k,\ell}\rvert \leq N$ such that
$$
\|\rho-\rho'- \langle n_{k,\ell},\alpha\rangle\|_{\R/\Z}< \epsilon^{\sigma}.
$$
By (\ref{g1}), for any $\rho\in E_{k}$ and $\rho'\in E_{\ell}$, we have
$$
\|\rho-\rho'- \langle n_{k,\ell},\alpha\rangle \|_{\R/\Z}<4m\epsilon^{\sigma}.
$$
Let  $\mathcal{R}_N\subset\left\{(k,\ell)|1\leq k<\ell\leq r\right\}$ be the set of all indexes such that $A_{k}$ and $A_{\ell}$ are $N$-resonant.

\begin{Lemma}
[Uniqueness of $N$-resonance]\label{unique} $n_{k,\ell}$ is the unique resonant site within length $\epsilon^{-\frac{\sigma}{2\tau}}\gg N$.
\end{Lemma}
\begin{pf}
Indeed, if there exists  $n'_{k,\ell}\neq n_{k,\ell}$ satisfying $\|\rho-\rho'- \langle n'_{k,\ell},\alpha\rangle\|_{\R/\Z}<4m\epsilon^{\sigma}$, then by the Diophantine condition of $\alpha$, we have
$$
\frac{\kappa}{| n'_{k,\ell}-n_{k,\ell}|^{\tau}}\leq \| \langle n'_{k,\ell}-n_{k,\ell},\alpha \rangle\|_{\R/\Z}<8m\epsilon^{\sigma},
$$
which implies that
$$\lvert n'_{k,\ell} \rvert>(8m)^{-\frac{1}{\tau}}\kappa^{\frac{1}{\tau}}\epsilon^{-\frac{\sigma}{\tau}}-\epsilon^{-\frac{\sigma}{2\tau}}\geq \epsilon^{-\frac{\sigma}{2\tau}}.$$
The last inequality holds since $\epsilon\leq2^{-\frac{2\tau}{\sigma}}(4m+1)^{-\frac{2}{\sigma}}\kappa^{\frac{2}{\sigma}}$.
\end{pf}

Now we are in the position to characterize  $\tilde{f}^{re}(x)$.
\begin{Lemma}[Structure of resonances]
There exists $1\leq i_0\leq m^2$ such that
 \begin{align*}
\tilde{f}^{re}(x)&=\tilde{f}^{re}_0+\tilde{f}^{re}_1(x)+\tilde{f}^{re}_2(x),
\end{align*}
\begin{align*}
\tilde{f}^{re}_0=\sum\limits_{\ell=1}^r\widehat{\tilde{f}^{re}_{\ell,\ell}}(0),\ \ \tilde{f}^{re}_2(x)=\sum_{\lvert n \rvert\geq N^{i_0+1}}\widehat{\tilde{f}^{re}}(n)e^{2\pi i\langle n,x\rangle},
\end{align*}
\begin{align*}
\tilde{f}^{re}_1(x)=\sum\limits_{(k,\ell)\in \mathcal{R}_{N^{i_0}}}\left(\widehat{\tilde{f}^{re}_{k,\ell}}(n_{k,\ell})e^{2\pi i\langle n_{k,\ell},x\rangle}+\widehat{\tilde{f}^{re}_{\ell,k}}(-n_{k,\ell})e^{-2\pi i\langle n_{k,\ell},x\rangle}\right).
\end{align*}
\end{Lemma}
\begin{pf}
Since $\Lambda_N\in \mathcal{B}_h^{nre}(\epsilon^{\frac{2}{5}})$, we have
\begin{align*}
\tilde{f}^{re}(x)=\sum\limits_{1\leq k,\ell\leq r}\sum_{n\in \Z^d\backslash\Lambda_{k,\ell}(\epsilon^{\sigma})}\widehat{\tilde{f}^{re}_{{k},\ell}}(n)e^{2\pi i\langle n,x\rangle}.
\end{align*}
Let $\mathcal{N}_i:=[N^i,N^{i+1})$ for $1\leq i\leq m^2$.  There are at most $m(m-1)$ different $n_{k,\ell}$'s, and $m^2$ different $\mathcal{N}_i$'s, thus by pigeonhole principle, there exists $1\leq i_0\leq m^2$ such that $n_{k,\ell}\notin \mathcal{N}_{i_0}$. I.e.,
$$
\| E_{k}-E_{\ell}-\langle n,\alpha\rangle\|_{\R/\Z} \geq \epsilon^{\sigma}, \ \  N^{i_0}\leq n<N^{i_0+1},\ \  \forall 1\leq k,\ell\leq r,
$$
which means there is no resonance in a large scale. We define
\begin{align}\label{eqneed1}
\tilde{f}^{re}_2(x)=\sum_{\lvert n \rvert\geq N^{i_0+1}}\widehat{\tilde{f}^{re}}(n)e^{2\pi i\langle n,x\rangle},
\end{align}
then it follows
\begin{align}\label{eqneed}
\tilde{f}^{re}(x)-\tilde{f}^{re}_2(x)=\sum\limits_{1\leq k,\ell\leq r}\sum_{n\in \{\Z^d\backslash\Lambda_{k,\ell}(\epsilon^{\sigma})\}\cap\{n|\lvert n \rvert\leq N^{i_0}\}}\widehat{\tilde{f}^{re}_{{k},\ell}}(n)e^{2\pi i\langle n,x\rangle}.
\end{align}
By the Diophantine condition on the frequency $\alpha$ and \eqref{g1}, for any $\rho,\rho'\in E_{\ell}$ ($1\leq \ell\leq r$) and $0<n\leq N^{i_0}$, we have
\begin{align}\label{est3}
\|\rho-\rho'- \langle n,\alpha\rangle\|_{\R/\Z} &\geq \frac{\kappa}{|n|^{\tau}}-4m\epsilon^{\sigma}\geq \kappa N^{-m^2\tau}-4m\epsilon^{\sigma}\geq \epsilon^{\sigma},
\end{align}
the last inequality holds since $\epsilon\leq \frac{2^{-\frac{2\tau}{\sigma}}(4m+1)^{-\frac{2}{\sigma}}\kappa^{\frac{2}{\sigma}}}{|A|^{\frac{500dm^3\tau}{\sigma}}}(h-h')^{\frac{500dm^3\tau}{\sigma}}$. It follows from \eqref{est3} that
\begin{equation}\label{inter1}
\left\{\Z^{d}\backslash\Lambda_{\ell,\ell}(\epsilon^{\sigma})\}\cap\{n\in \Z^{d}:\lvert n \rvert\leq N^{i_0}\right\}=\{0\},\ \ \ell=1,\cdots,r.
\end{equation}
On the other hand, by Lemma \ref{unique}, we have
\begin{equation}\label{inter2}
\left\{\Z^{d}\backslash\Lambda_{k,\ell}(\epsilon^{\sigma})\}\cap\{n\in \Z^{d}:\lvert n \rvert\leq N^{i_0}\right\}=\{\pm n_{k,\ell}\}, \ \  k\neq  \ell,
\end{equation}
since $N^{m^2}\leq \epsilon^{-\frac{\sigma}{2\tau}}$ by \eqref{initial_assump}.

\eqref{eqneed1}, \eqref{eqneed},\eqref{inter1} and \eqref{inter2} finish the proof.
\end{pf}  

\bigskip
\noindent 
{\bf Step 4. Eliminating the  lower order resonant terms}.
\begin{Lemma}
There exists a family $m_1,\cdots, m_{r}$ with $\max\limits_{1\leq\ell\leq r}|m_\ell|\leq rN^{i_0}$ such that
\begin{flalign*}
m_k-m_\ell=-n_{k,\ell},\ \  (k,\ell)\in \mathcal{R}_{N^{i_0}}.
\end{flalign*}
\end{Lemma}
\begin{pf}
We prove this by induction, assume for $\mathcal{R}_{N^{i_0}}\cap\{(k,\ell)|1\leq k<\ell\leq r-1\}$, the above lemma holds, which means there exists a family $m_1\cdots m_{r-1}$ with $\max\limits_{1\leq\ell\leq r}|m_\ell|\leq (r-1)N^{i_0}$ such that
\begin{flalign*}
m_k-m_\ell=-n_{k,\ell}, \ \ (k,\ell)\in\mathcal{R}_{N^{i_0}}\cap \{(k,\ell)|1\leq k<\ell\leq r-1\}.
\end{flalign*}
We consider $\mathcal{R}_{N^{i_0}}$. There are two possible cases.\\
Case I: There exists $1\leq \ell\leq r-1$, such that $(\ell,r)\in\mathcal{R}_{N^{i_0}}$. Thus there exists $0<|n_{\ell,r}|\leq N^{i_0}$ such that
$$
\|\rho-\rho'-\langle n_{\ell,r},\alpha\rangle\|_{\R/\Z}\leq \epsilon^{\sigma},
$$
for some $\rho\in E_{\ell}$ and $\rho'\in E_{{r}}$. On the other hand, if there exists $1\leq \ell'\leq r-1$, such that $(\ell',r)\in\mathcal{R}_{N^{i_0}}$, then there exists $0<|n_{\ell',r}|\leq N^{i_0}$ such that
$$
\|\tilde{\rho}-\tilde{\rho}'-\langle n_{\ell',r},\alpha\rangle\|_{\R/\Z}\leq \epsilon^{\sigma},
$$
for some $\tilde{\rho}\in E_{{\ell'}}$ and $\tilde{\rho}'\in E_{{r}}$. This implies that
$$
\|\rho-\tilde{\rho}-\langle n_{\ell,r}-n_{\ell',r},\alpha\rangle\|_{\R/\Z}\leq 4m\epsilon^{\sigma}.
$$
Thus
$$
n_{\ell,r}-n_{\ell',r}=n_{\ell,\ell'}.
$$
Let $m_{r}=m_\ell+n_{\ell,r}$, then
$$
m_{\ell'}-m_{r}=m_{\ell'}-m_\ell-n_{\ell,r}=n_{\ell,\ell'}-n_{\ell,r}=-n_{\ell',r},
$$
for $\ell'\neq r$, with the estimate
$$
|m_{r}|\leq |m_\ell|+|n_{\ell,r}|\leq (r-1)N^{i_0}+N^{i_0}=rN^{i_0}.
$$
\\
Case II: $(\ell,r)\notin\mathcal{R}_{N^{i_0}}$ for any $1\leq \ell\leq r-1$, let $m_{r}=0$, we get the result.
\end{pf}

Define the $\Z^d$-periodic rotation $Q(x)$ as below:
$$
Q(x)=diag\{e^{-2\pi i\langle m_1,x\rangle}I_{n_1},\cdots,e^{-2\pi i\langle m_r,x\rangle}I_{n_r}\}.
$$
So we have \begin{equation}\label{esti-Q}\lvert Q\rvert_{h_+}\leq e^{rN^{i_0}h_+}\leq  e^{mN^{i_0}h_+}.\end{equation}
One can also show that
\begin{equation}\label{star}
Q(x+\alpha)\widetilde{A}e^{\tilde{f}^{re}(x)}Q^{-1}(x)=\widehat{A}e^{\hat{f}(x)},
\end{equation}
where
$$
\widehat{A}=Q(x+\alpha)diag\{A_{1},\cdots, A_{r}\}Q^{-1}(x)=diag\{A_{1}e^{-2\pi i\langle m_1,\alpha\rangle},\cdots,A_{r}e^{-2\pi i\langle m_r,\alpha\rangle}\},
$$
and
$$\hat{f}(x)=Q(x)\tilde{f}^{re}(x)Q^{-1}(x)=Q(x)\tilde{f}_0^{re}Q^{-1}(x)+Q(x)\tilde{f}^{re}_1(x)Q^{-1}(x)+Q(x)\tilde{f}^{re}_2(x)Q^{-1}(x).$$
Moreover,
\begin{flalign*}
Q(x)\tilde{f}_0^{re}Q^{-1}(x) &=\tilde{f}_0^{re} \in gl(m,\C),
\end{flalign*}
\begin{flalign*}
Q(x)\tilde{f}^{re}_1(x)Q^{-1}(x)=\sum\limits_{(k,\ell)\in \mathcal{R}_{N^{i_0}}}\left(\widehat{\tilde{f}^{re}_{k,\ell}}(n_{k,\ell})+\widehat{\tilde{f}^{re}_{\ell,k}}(-n_{k,\ell})\right)\in gl(m,\C).
\end{flalign*}
Denote
\begin{flalign*}
L=Q(x)\tilde{f}_0^{re}Q^{-1}(x)+Q(x)f^{re}_1(x)Q^{-1}(x),\ \ F(x)=Q(x)f^{re}_2(x)Q^{-1}(x), \ \ \bar{B}(x)=Q(x)e^{Y(x)}P.
\end{flalign*}
By \eqref{new1}, \eqref{new2} and \eqref{star}, we have
\begin{equation}\label{con1}
\bar{B}(x+\alpha)Ae^{f(x)}\bar{B}^{-1}(x)=\widehat{A}e^{\hat{f}(x)}=\widehat{A}e^{L+F(x)}.
\end{equation}
Moreover, we have the following estimates:
\begin{equation}\label{need_1}
|\bar{B}|_0 \leq |e^{Y}|_{h_+} |P| \leq 2C_0(2|A|\epsilon^{-\sigma})^{m(m+1)}\leq \epsilon^{-50m^2\sigma},
\end{equation}
\begin{equation*}
\label{F}|F|_{h_+} \leq \lvert Qf^{re}_2Q^{-1}\rvert_{h_+} \leq 2\epsilon^{\frac{9}{10}}e^{-N^{i_0+1}(h-h_+)}e^{2mN^{i_0}h_+}\leq 2\epsilon^{\frac{9}{10}}e^{-N^{i_0}(N(h-h_+)-2mh_+)}\leq \epsilon^3.
\end{equation*}
\begin{equation}\label{star1}
|L|\leq |\tilde{f}^{re}_0|+\sum\limits_{(k,\ell)\in \mathcal{R}_{N^{i_0}}}\left(|\widehat{\tilde{f}^{re}_{k,\ell}}(n_{k,\ell})|+|\widehat{\tilde{f}^{re}_{\ell,k}}(-n_{k,\ell})|\right)\leq 2m^2|\tilde{f}^{re}|_h\leq 2m^2\epsilon^{1-200m^2\sigma}.
\end{equation}
Direct computation shows that
\begin{equation}\label{impl}
e^{L+F(x)}=e^L+\mathcal{O}(F(x))=e^L(Id+e^{-L}\mathcal{O}(F(x)))=e^L e^{f_+{(x)}}.
\end{equation}
It immediately implies that
$$
\lvert f_+\rvert_{h'}\leq 2 |F|_{h_+}\leq \epsilon_+.
$$
Thus we can rewrite  $(\ref{con1})$ as
\begin{equation}\label{eq_2}
\bar{B}(x+\alpha)(Ae^{f(x)})\bar{B}^{-1}(x)=A_{+}e^{f_+(x)},
\end{equation}
with
\begin{equation}\label{eigen_A}
A_+=diag\{A_{1}e^{-2\pi i\langle m_1,\alpha\rangle},\cdots,A_{r}e^{-2\pi i\langle m_r,\alpha\rangle}\}e^L.
\end{equation}
\

\ 

Combining the above four steps, let $B_{j+1}(x)=\bar{B}(x)B_j(x)$, $A_{j+1}=A_+$ and $f_{j+1}=f_+$, then $B_{j+1}\in C^\omega_{h_{j+1}}(\T^d,GL(m,\C))$, by \eqref{eq_1} and \eqref{eq_2}, we have
\begin{equation}\label{star2}
B_{j+1}(x+\alpha)A_0e^{f_0(x)}B^{-1}_{j+1}(x)=A_{j+1}e^{f_{j+1}(x)}.
\end{equation}
By \eqref{ess3'} and \eqref{need_1}, we have
\begin{equation}\label{star3}
|B_{j+1}|_{0}\leq |B_j|_0|\bar{B}|_0\leq\epsilon_{j-1}^{-200m^2\sigma}\epsilon_j^{-50m^2\sigma}\leq \epsilon_j^{-200m^2\sigma}.
\end{equation}
By \eqref{eigen_A}, Proposition \ref{eigen_per} and \eqref{star1}, we can permute the eigenvalues of $A_{j+1}$ as $\{e^{-2\pi i\rho^{j+1}_\ell}\}_{\ell=1}^m$, such that
\begin{equation}\label{star4}
\big|\Im\rho_\ell^{j+1}-\Im\rho^j_\ell\big|\leq C(m,A_j)(2m^2\epsilon_j^{(1-200m^2\sigma)})^{\frac{1}{m}}\leq C(m,A_0)\epsilon_j^{\frac{1}{m}-200m^2\sigma}, \ \ 1\leq \ell \leq m.
\end{equation}
On the other hand, it is obvious that
\begin{equation}\label{star5}
\left|\sum\limits_{\ell=1}^m\Im\rho_\ell^{j+1}-\sum\limits_{\ell=1}^m\Im\rho^j_\ell\right|\leq C(m,A_0)\epsilon_j^{\frac{1}{2}}.
\end{equation}
\eqref{star2}, \eqref{star3}, \eqref{star4} and \eqref{star5} finish the proof.\qed
\begin{Corollary}\label{coro_1}
Assume $\left(\alpha , A\right)$ is $C_{h}^\omega$ almost reducible, for any given $0<\sigma<\frac{1}{500m^3}$ and $0<h'<h$,  there exist $B_{j}\in C^{\omega}_{h'}(\T^{d}$,
$GL(m,\C))$, $A_{j}\in GL(m,\C)$ and $f_{j}\in C^{\omega}_{h'}(\T^d,gl(m,\C))$, such that
$$
B_{j}(x+\alpha)A(x)B^{-1}_{j}(x)=A_{j}e^{f_j(x)}.
$$
Moreover, we have the following estimates
\begin{equation*}
|f_{j}|_{h'}\leq \epsilon_{j},\ \ |B_{j}|_0\leq C(\alpha,A)\epsilon_{j-1}^{-200m^2\sigma},
\end{equation*}
\begin{align*}
\left|\sum\limits_{\ell=1}^m\Im\rho_{\ell}^j-\sum\limits_{\ell=1}^m\Im\rho_{\ell}^{j-1}\right|\leq C(m,A)\epsilon_{j-1}^{\frac{1}{2}},
\end{align*}
\begin{equation*}
\left| \Im\rho_{\ell}^j-\Im\rho_{\ell}^{j-1}\right|\leq C(m,A)\epsilon^{\frac{1}{m}-200m^2\sigma}_{j-1},\ \ 1\leq \ell\leq m,
\end{equation*}
where $\left\{e^{-2\pi i\rho_\ell^j}\right\}_{\ell=1}^m$ are the eigenvalues of $A_j$ and $\epsilon_j= \epsilon_0^{2^j}$.
\end{Corollary}
\begin{pf}
By the assumption, $(\alpha,A)$ is $C^\omega_h$ almost reducible, there exist $B\in C^{\omega}_{h}(\T^{d}$,
$GL(m,\C))$, $A_0\in GL(m,\C)$, $f_0\in C^{\omega}_{h}(\T^d,gl(m,\C))$, such that
$$ 
B(x+\alpha)A(x)B^{-1}(x)=A_0e^{f_0(x)},
$$
with $|f_0|_h\leq\epsilon_{0}\leq \epsilon(\alpha,h,h',m,|A|,\sigma)$ which is defined in Theorem \ref{thm3}.  Applying Theorem \ref{thm3} to $(\alpha, A_0e^{f_0(\cdot)})$,  there exist $\bar{B}_{j}\in C^{\omega}_{h'}(\T^{d}$,
$GL(m,\C))$, $A_{j}\in GL(m,\C)$, $f_{j}\in C^{\omega}_{h'}(\T^d, gl(m,\C))$, such that
$$
\bar{B}_{j}(x+\alpha)(A_0e^{f_{0}(x)})\bar{B}^{-1}_{j}(x)=A_{j}e^{f_j(x)}.
$$
Moreover, we have the following estimates
\begin{equation*}
|f_{j}|_{h'}\leq \epsilon_{j},\ \ |\bar{B}_{j}|_0\leq \epsilon_{j-1}^{-200m^2\sigma},
\end{equation*}
\begin{align*}
\left|\sum\limits_{\ell=1}^m\Im\rho_{\ell}^j-\sum\limits_{\ell=1}^m\Im\rho_{\ell}^{j-1}\right|\leq C(m,A_0)\epsilon_{j-1}^{\frac{1}{2}},
\end{align*}
\begin{equation*}
\left| \Im\rho_{\ell}^j-\Im\rho_{\ell}^{j-1}\right|\leq C(m,A_0)\epsilon^{\frac{1}{m}-200m^2\sigma}_{j-1},\ \ 1\leq \ell\leq m,
\end{equation*}
where $\left\{e^{-2\pi i\rho_\ell^j}\right\}_{\ell=1}^m$ are the eigenvalues of $A_j$ and $\epsilon_j= \epsilon_0^{2^j}$.

The desired results follow if we denote $B_j(x)=\bar{B}_j(x)B(x)$. 
\end{pf}

\section{Proof of the main Theorems}
In this section, we give the proof of the main theorems. Our main ideas are based on the following basic facts.
\begin{Proposition}\label{final1}
Assume that $(\alpha,A)/(\alpha,\widetilde{A})$ is $C^\omega_h$ almost reducible to $(\alpha,A_\infty)/(\alpha,\widetilde{A}_\infty)$. Denote the eigenvalues of $A_\infty/\widetilde{A}_\infty$ by $\{e^{-2\pi i\rho^\infty_j}\}_{j=1}^m/\{e^{-2\pi i\tilde{\rho}^\infty_j}\}_{j=1}^m$ respectively. If
\begin{equation}\label{eigenle}
\lvert \Im\rho^\infty_{j}-\Im\tilde{\rho}^\infty_{j}\lvert\leq \epsilon,\ \ 1\leq j\leq m,
\end{equation}
then
\begin{equation}\label{le1}
|L_j(\alpha,A)-L_j(\alpha,\widetilde{A})|\leq 100\epsilon,\ \ 1\leq j\leq m.
\end{equation}
\end{Proposition}
\begin{pf}
We prove \eqref{le1} inductively based on \eqref{eigenle}. For $m=1$, by Proposition \ref{leconstant}, it is obvious. Now we assume \eqref{le1} holds for $m\leq k$. We consider the case $m=k+1$, again by Proposition \ref{leconstant}, we can assume
$$
L_1(\alpha,A)=\max_j2\pi\Im \rho^\infty_j=2\pi\Im \rho^\infty_{j_0},  \ \  L_1(\alpha,\widetilde{A})=\max_j2\pi\Im \tilde{\rho}^\infty_j=2\pi\Im \tilde{\rho}^\infty_{j_1}.
$$
We distinguish into two cases,\\
Case I: $j_0=j_1$, by \eqref{eigenle},
$$
|L_1(\alpha,A)-L_1(\alpha,\widetilde{A})|=2\pi|\Im \rho^\infty_{j_0}-\Im \tilde{\rho}^\infty_{j_1}|\leq 2\pi\epsilon.
$$
By induction, we have
\begin{equation*}
|L_j(\alpha,A)-L_j(\alpha,\widetilde{A})|\leq 100\epsilon,\ \ 2\leq j\leq k+1.
\end{equation*}

\noindent
Case II: $j_0\neq j_1$, by \eqref{eigenle},
$$
|\Im \rho^\infty_{j_0}-\Im \tilde{\rho}^\infty_{j_0}|\leq \epsilon,\  \ |\Im \rho^\infty_{j_1}-\Im \tilde{\rho}^\infty_{j_1}|\leq \epsilon.
$$
On the other hand,
$$
L_1(\alpha,A)=\max_j2\pi\Im \rho^\infty_j=2\pi\Im \rho^\infty_{j_0},  \ \  L_1(\alpha,\widetilde{A})=\max_j2\pi\Im \tilde{\rho}^\infty_j=2\pi\Im \tilde{\rho}^\infty_{j_1}.
$$
We must have 
$$
|\Im \rho^\infty_{j_0}-\Im \rho^\infty_{j_1}|\leq 50\epsilon.
$$
Thus
$$
|L_1(\alpha,A)-L_1(\alpha,\widetilde{A})||=|\Im \rho^\infty_{j_0}-\Im \tilde{\rho}^\infty_{j_1}|\leq|\Im \rho^\infty_{j_0}-\Im \rho^\infty_{j_1}|+ |\Im \rho^\infty_{j_1}-\Im \tilde{\rho}^\infty_{j_1}|\leq 100\epsilon.
$$
By induction, we have
\begin{equation*}
|L_j(\alpha,A)-L_j(\alpha,\widetilde{A})|\leq 100\epsilon,\ \ 2\leq j\leq k+1.
\end{equation*}
Thus we finish the proof.
\end{pf}
For any given $A_0\in GL(m,\C)$ with eigenvalues $\{e^{-2\pi i\rho_\ell^0}\}_{\ell=1}^m$, we assume that $\{\rho^0_\ell\}_{1\leq \ell\leq m}$ can be grouped as  $\bigcup\limits_{k=1}^r E_{k}$ with $\# E_{k}=n_k$, such that
\begin{equation}\label{c1}
|\rho-\rho'|\geq \delta>0,\ \ \rho\in E_{k},\rho'\in E_{\ell}, \ \ \forall k\neq \ell.
\end{equation}
Then we have the following refinements of Theorem \ref{thm3}.
\begin{Proposition}\label{final2}
Assume $\alpha\in {\rm DC}_d(\kappa,\tau)$, $0<\sigma<\frac{1}{500m^3}$ and $0<h'<h$, there exists $\epsilon_0=\epsilon(\alpha,h,h',m, |A_0|,\sigma,\delta)$ such that if $|f_0|_h\leq \epsilon_0$, then $(\alpha , A_0e^{f_0(\cdot)})$ is $C_{h'}^\omega$ almost reducible to $(\alpha,A_\infty)$. Moreover, $\{e^{-2\pi i\rho^\infty_\ell}\}_{\ell=1}^m$, the eigenvalues of $A_\infty$ satisfy
\begin{equation}\label{eigenle1}
\left\lvert \Im\rho^\infty_{\ell}-\Im\rho^0_{\ell}\right\lvert\leq C(m,\delta,A_0)\epsilon_0^{\frac{1}{n_k}-200m^2\sigma},\ \ \sum\limits_{j=1}^{k-1}n_j+1\leq \ell\leq \sum\limits_{j=1}^{k}n_j,\ \ 1\leq k\leq r.
\end{equation}
\begin{equation}\label{eigenle2}
\left\lvert \sum\limits_{\ell=1}^{n_k}\Im\rho^\infty_{\ell}-\sum\limits_{\ell=1}^{n_k}\Im\rho^0_{\ell}\right\lvert\leq C(m,\delta,A_0)\epsilon_0^{\frac{1}{2}},\ \ 1\leq k\leq r.
\end{equation}
\end{Proposition}
\begin{pf}
By Lemma $A^1$ in \cite{eliasson}, there exists $P\in GL(m,\mathbb{C})$ such that
\begin{equation*}
PA_0P^{-1}=\widetilde{A}_0:=diag\{A_0^{1},A_0^{2},\cdots, A_0^{r}\},
\end{equation*}
with $\|P\|\leq C(\delta,m,|A_0|)$ and $spec(A_0^{k})=\{e^{-2\pi i\rho}|\rho\in E_{k}\}  (k=1,\cdots,r)$. 

Thus we can conjugate $\left(\alpha, A_0e^{f_0(\cdot)}\right)$ to a new cocycle,
\begin{equation*}
PA_0e^{f_0(x)}P^{-1}=diag\{A_0^{1}\cdots A_0^{r}\}e^{Pf_0(x)P^{-1}}:=\widetilde{A}_0e^{\tilde{f}_0(x)},
\end{equation*}
with
\begin{equation*}
|\tilde{f}|_h \leq |P||f_0|_h|P^{-1}|\leq C(m,\delta,|A_0|)|f_0|_h.
\end{equation*}
Now we define
\begin{align*}
\Lambda=\left\{f\in C^{\omega}_{h}(\T^{d},gl(m,\C))\mid f(x)=\sum\limits_{1\leq k\neq \ell\leq r}\hat{f}_{k,\ell}(n)e^{2\pi i\langle n,x\rangle}\right\},
\end{align*}    
where $\hat{f}_{k,\ell}(n)=P_{[\sum_{j=1}^{k-1}n_{j}+1,\sum_{j=1}^{k}n_j]}\hat{f}(n)P_{[\sum_{j=1}^{\ell-1}n_{j}+1,\sum_{j=1}^\ell n_j]}$ and $P_I:[1,m] \rightarrow I$ is a projection. Similar to Lemma \ref{non-eli}, one can verify that for any $Y\in\Lambda$,
$$
\lvert \widetilde{A}_0^{-1}Y(\cdot+\alpha)\widetilde{A}_0-Y(\cdot)\rvert_h\geq c(m)\delta\left(1+m\left(|\widetilde{A}_0|+|\widetilde{A}_0|\right)\delta^{-1}\right)^{-2m}\lvert Y\rvert_h\geq  c(m,A_0,\delta)\lvert Y\rvert_h.
$$
Thus $\Lambda\subset \mathcal{B}_h^{nre}(c)$. Choose $\epsilon_0$ sufficiently small such that $c\geq 13|\widetilde{A}_0|^2C^2\epsilon_0^2$, by Lemma \ref{lem2}, there exists $B\in C^\omega(\T^d,GL(m,\mathbb{C}))$ such that
\begin{equation*}
B(x+\alpha)\widetilde{A}_0e^{\tilde{f}_0(x)}B^{-1}(x)=\widetilde{A}_0e^{\tilde{f}^{re}_0(x)}:=diag\left\{A_0^{1}e^{f_0^{1}(x)},A_0^{2}e^{f_0^{2}(x)},\cdots,A_0^{r}e^{f_0^{r}(x)}\right\},
\end{equation*}
with estimates
\begin{equation}
|f_0^{k}|_{h}\leq C(m,A_0,\delta)|f_0|_h\leq C(m,A_0,\delta)\epsilon_0,\ \  k=1,2,\cdots,r, \ \ |B|_0\leq 2.
\end{equation}
Assume $\epsilon_0\leq \epsilon(\alpha,h,h',m, |A_0|,\sigma)\delta^{100}$ where $\epsilon(\alpha,h,h',m, |A_0|,\sigma)$ is defined in Theorem \ref{thm3}. Applying Theorem \ref{thm3} to $(\alpha,A_0^{k}e^{f_0^{k}(x)})$ $(k=1,\cdots,r)$,  there exist $B^{k}_{j}\in C^{\omega}_{h'}(\T^{d}$,
$GL(n_k,\C))$, $A^{k}_{j}\in GL(n_k,\C)$, $f^{k}_{j}\in C^{\omega}_{h'}(\T^d,gl(n_k,\C))$, such that
$$
B^{k}_{j}(x+\alpha)(A^{k}_0e^{f_0^{k}(x)})(B^{k}_{j})^{-1}(x)=A^{k}_{j}e^{f^{k}_j(x)}.
$$
Moreover, we have the following estimates
\begin{equation*}
\left| \sum\limits_{\ell=\sum\limits_{j=1}^{k-1}n_j+1}^{\sum\limits_{j=1}^{k}n_j}\Im\rho_{\ell}^j- \sum\limits_{\ell=\sum\limits_{j=1}^{k-1}n_j+1}^{\sum\limits_{j=1}^{k}n_j}\Im\rho_{\ell}^{j-1}\right|\leq C(m,A_0)\epsilon^{\frac{1}{2}}_{j-1},
\end{equation*}
\begin{equation*}
\left| \Im\rho_{\ell}^j-\Im\rho_{\ell}^{j-1}\right|\leq C(m,A_0)\epsilon^{\frac{1}{n_k}-200m^2\sigma}_{j-1},\ \ \sum\limits_{j=1}^{k-1}n_j+1\leq \ell\leq \sum\limits_{j=1}^{k}n_j,
\end{equation*}
where $\left\{e^{-2\pi i\rho_\ell^j}\right\}_{\ell=\sum\limits_{j=1}^{k-1}n_j+1}^{\sum\limits_{j=1}^{k}n_j}$ are the eigenvalues of $A^{k}_j$ and $\epsilon_j= \epsilon_0^{2^j}$.

It follows that
\begin{equation*}
\left| \sum\limits_{\ell=\sum\limits_{j=1}^{k-1}n_j+1}^{\sum\limits_{j=1}^{k}n_j}\Im\rho_{\ell}^\infty- \sum\limits_{\ell=\sum\limits_{j=1}^{k-1}n_j+1}^{\sum\limits_{j=1}^{k}n_j}\Im\rho_{\ell}^{0}\right|\leq \sum\limits_{j=1}^\infty\left| \sum\limits_{\ell=\sum\limits_{j=1}^{k-1}n_j+1}^{\sum\limits_{j=1}^{k}n_j}\Im\rho_{\ell}^j-\sum\limits_{\ell=\sum\limits_{j=1}^{k-1}n_j+1}^{\sum\limits_{j=1}^{k}n_j}\Im\rho_{\ell}^{j-1}\right|\leq C(m,\delta,A_0)\epsilon_0^{\frac{1}{2}},
\end{equation*}
\begin{align*}
\left\lvert \Im\rho^\infty_{\ell}-\Im\rho^0_{\ell}\right\lvert&\leq \sum\limits_{j=1}^\infty\left| \Im\rho_{\ell}^j-\Im\rho_{\ell}^{j-1}\right|\leq \sum\limits_{j=1}^\infty C(m,\delta,A_0)\epsilon^{\frac{1}{n_k}-200m^2\sigma}_{j-1}\\
&\leq C(m,\delta,A_0)\epsilon_0^{\frac{1}{n_k}-200m^2\sigma},\ \ \sum\limits_{j=1}^{k-1}n_j+1\leq \ell\leq \sum\limits_{j=1}^{k}n_j.
\end{align*}
Thus we finish the whole proof.
\end{pf}
If we further assume the eigenvalues of $A_0$ satisfy
\begin{equation}\label{c2}
\rho-\rho'\geq \delta>0,\ \ \rho\in E_{k},\rho'\in E_{\ell}, \ \ \forall k< \ell.
\end{equation}
We immediately have the following corollary.
\begin{Corollary}\label{final3}
Assume $\alpha\in {\rm DC}_d(\kappa,\tau)$ and $0<\sigma<\frac{1}{500m^3}$. There exists $\epsilon_0=\epsilon(\alpha,h,m, |A_0|,\sigma,\delta)$ such that if $|f_0|_h\leq \epsilon_0$, then
\begin{equation*}
\left\lvert L_\ell(\alpha, A_0e^{f_0(\cdot)})-\Im\rho^0_{\ell}\right\lvert\leq C(m,\delta,A_0)\epsilon_0^{\frac{1}{n_k}-200m^2\sigma},\ \ \sum\limits_{j=1}^{k-1}n_j+1\leq \ell\leq \sum\limits_{j=1}^{k}n_j,\ \ 1\leq k\leq r.
\end{equation*}
\begin{equation*}
\left\lvert \sum\limits_{\ell=1}^{n_k} L_{\ell}(\alpha,A_0e^{f_0(\cdot)})-\sum\limits_{\ell=1}^{n_k}\Im\rho^0_{\ell}\right\lvert\leq C(m,\delta,A_0)\epsilon_0^{\frac{1}{2}},\ \ 1\leq k\leq r.
\end{equation*}
\end{Corollary}
\begin{pf}
Take $h'=\frac{h}{2}$ and $\epsilon_0=\epsilon(\alpha,h,h/2,m, |A_0|,\sigma,\delta)$ which is defined in Proposition \ref{final2},  then $(\alpha , A_0e^{f_0(\cdot)})$ is $C_{h'}^\omega$ almost reducible to $(\alpha,A_\infty)$. Moreover, $\{e^{-2\pi i\rho^\infty_\ell}\}_{\ell=1}^m$, the eigenvalues of $A_\infty$ satisfy
\begin{equation}\label{eigenle1'}
\left\lvert \Im\rho^\infty_{\ell}-\Im\rho^0_{\ell}\right\lvert\leq C(m,\delta,A_0)\epsilon_0^{\frac{1}{n_k}-200m^2\sigma},\ \ \sum\limits_{j=1}^{k-1}n_j+1\leq \ell\leq \sum\limits_{j=1}^{k}n_j,\ \ 1\leq k\leq r.
\end{equation}
\begin{equation}\label{eigenle2'}
\left\lvert \sum\limits_{\ell=1}^{n_k}\Im\rho^\infty_{\ell}-\sum\limits_{\ell=1}^{n_k}\Im\rho^0_{\ell}\right\lvert\leq C(m,\delta,A_0)\epsilon_0^{\frac{1}{2}},\ \ 1\leq k\leq r.
\end{equation}
By \eqref{c2} and \eqref{eigenle1'},  $\{\rho^\infty_\ell\}_{\ell=1}^m$ can be grouped as  $\bigcup\limits_{k=1}^r E^\infty_{k}$ with $\# E^\infty_{k}=n_k$, such that
\begin{equation}\label{c2'}
\rho-\rho'\geq \frac{\delta}{2},\ \ \rho\in E^\infty_{k},\rho'\in E^\infty_{\ell}, \ \ \forall k< \ell.
\end{equation}
By \eqref{c2'} and Proposition \ref{leconstant}, we have
\begin{equation}\label{c2''}
\left\{L_\ell(\alpha,A_0e^{f_0(\cdot)})\right\}_{\ell=\sum\limits_{j=1}^{k-1}n_j+1}^{\sum\limits_{j=1}^{k}n_j}=E^\infty_k, \ \ 1\leq k\leq r.
\end{equation}
The desired result follows from \eqref{eigenle1'}, \eqref{eigenle2'}, \eqref{c2''} and Proposition \ref{final1}. 
\end{pf}
\subsection{Proof of Theorem \ref{thm1}} We assume $(\alpha,A)/(\alpha,\widetilde{A})$ is almost reducible to $(\alpha,A_\infty)/(\alpha,\widetilde{A}_\infty)$ respectively. We divide the proof into the following three steps.

\noindent{\bf Step 1. Group the eigenvalues.}
Since $(\alpha,A)$ is $C^\omega_h$ almost reducible, by Corollary \ref{coro_1}, for any $\e>0$, there exist $B_{j}\in C^{\omega}_{h'}(\T^{d}$,
$GL(m,\C))$, $A_{j}\in GL(m,\C)$ and $f_{j}\in C^{\omega}_{h'}(\T^d,gl(m,\C))$, such that
$$
B_{j}(x+\alpha)A(x)B^{-1}_{j}(x)=A_{j}e^{f_j(x)}.
$$
Moreover, we have the following estimates
\begin{equation}\label{n1}
|f_{j}|_{h'}\leq \epsilon_{j},\ \ |B_{j}|_0\leq C(\alpha,A)\epsilon_{j-1}^{-200m^2\e},
\end{equation}
\begin{equation}\label{n3}
\left| \Im\rho_{\ell}^j-\Im\rho_{\ell}^{j-1}\right|\leq C(m,A)\epsilon^{\frac{1}{m}-200m^2\e}_{j-1},\ \ 1\leq \ell\leq m,
\end{equation}
where $\left\{e^{-2\pi i\rho_\ell^j}\right\}_{\ell=1}^m$ are the eigenvalues of $A_j$ and $\epsilon_j= \epsilon_0^{2^j}$.

We assume there are $r$ distinct Lyapunov exponents of $(\alpha,A)$. By Proposition \ref{leconstant}, there are $r$ distinct $\Im \rho^\infty_\ell$'s. Thus $\{\rho^\infty_\ell\}_{\ell=1}^m$ can be grouped as $\bigcup\limits_{k=1}^r E^{\infty}_{k}$ with $\# E^{\infty}_{k}=n_k$, such that
\begin{equation*}
\rho-\rho'\geq 2\eta>0,\ \ \rho\in E^{\infty}_{k},\rho'\in E^{\infty}_{\ell}, \ \ k<\ell.
\end{equation*}
By \eqref{n3}, for $j$ sufficiently large depending on $\alpha,A,\eta$, we can group $\{\rho^{j+1}_\ell\}_{\ell=1}^m$ as $\bigcup\limits_{k=1}^r E^{j+1}_{k}$ with $\# E^{j+1}_{k}=n_k$, such that
\begin{equation*}
\rho-\rho'\geq \eta,\ \ \rho\in E^{j+1}_{k},\rho'\in E^{j+1}_{\ell}, \ \ k<\ell.
\end{equation*}
\noindent{\bf Step 2. Compare $(\alpha,A)$ and $(\alpha,\widetilde{A})$.} Denote $\delta=|B-A|_h$, we only need to consider the case that $\delta$ is sufficiently small. Assume that $\epsilon_0\leq \epsilon(\alpha,h,h',|A|,m,\e,\eta)$ which is defined in Proposition \ref{final2} and $\epsilon_{j+1}\leq \delta\leq \epsilon_j$. It is obvious that
\begin{align*}
B_{j+1}(x+\alpha)\widetilde{A}(x)B_{j+1}^{-1}(x)&=B_{j+1}(x+\alpha)(\widetilde{A}(x)-A(x))B_{j+1}^{-1}(x)+B_{j+1}(x+\alpha)A(x)B_{j+1}^{-1}(x)\\
&=B_{j+1}(x+\alpha)(B(x)-A(x))B_{j+1}^{-1}(x)+A_{j+1}e^{f_{j+1}(x)}\\
&:=A_{j+1}e^{\tilde{f}_{j+1}(x)}.
\end{align*}
By \eqref{n1}, we have
$$
|\tilde{f}_{j+1}|_h\leq 100\delta|B_{j+1}|_0^2+\epsilon_{j+1}\leq C(\alpha,A)\delta\epsilon_j^{-800m^2\e}\leq  C(\alpha,A)\delta^{1-800m^2\e}.
$$ 
\noindent{\bf Step 3. H\"older continuity of the Lyapunov exponents.} 
By Corollary \ref{final3}, we have 
\begin{enumerate}
\item $\left\lvert L_\ell\left(\alpha,A_{j+1}e^{f_{j+1}(\cdot)}\right)-\Im\rho^{j+1}_{\ell}\right\lvert\leq C(m,\delta,A)\delta^{\frac{1}{n_k}-200m^2\e},\ \ \sum\limits_{j=1}^{k-1}n_j+1\leq \ell\leq \sum\limits_{j=1}^{k}n_j,\ \ 1\leq k\leq r$.
\item $\left\lvert L_\ell\left(\alpha,A_{j+1}e^{\tilde{f}_{j+1}(\cdot)}\right)-\Im\rho^{j+1}_{\ell}\right\lvert\leq C(m,\delta,A)\delta^{\frac{1}{n_k}-800m^2\e},\ \ \sum\limits_{i=1}^{k-1}n_j+1\leq \ell\leq \sum\limits_{j=1}^{k}n_j,\ \ 1\leq k\leq r$.
\end{enumerate}
Thus 
$$
\left\lvert L_\ell\left(\alpha,A_{j+1}e^{f_{j+1}(\cdot)}\right)-L_\ell\left(\alpha,A_{j+1}e^{\tilde{f}_{j+1}(\cdot)}\right)\right\lvert\leq C(m,\delta,A)\delta^{\frac{1}{n_k}-800m^2\e},\ \ \sum\limits_{j=1}^{k-1}n_j+1\leq \ell\leq \sum\limits_{j=1}^{k}n_j,\ \ 1\leq k\leq r.
$$
By Proposition \ref{invariance}, we have 
\begin{align*}
|L_\ell(\alpha,A)-L_\ell(\alpha,\widetilde{A})|\leq C\delta^{\frac{1}{n_k}(1-800m^2\e)}, \ \ \sum\limits_{j=1}^{k-1}n_j+1\leq \ell\leq \sum\limits_{j=1}^{k}n_j, \ \ 1\leq k\leq r.
\end{align*}

\subsection{Proof of Theorem \ref{main}} We only need to consider the case that $E$ is located in a bounded set $\mathcal{B}$, otherwise $(\alpha,L_{E,W}^{\lambda^{-1}V})$ is uniformly hyperbolic and the IDS is smooth.  For any $E\in\mathcal{B}$, assume $\lambda^{-1}\leq \lambda_0(m,V,W)$, we can rewrite $(\alpha,L_{E,W}^{\lambda^{-1}V})$ as $(\alpha,L_{E,W}^{0}e^{f(x)})$ with $|f|_h\leq \lambda^{-\frac{1}{2}}$. We divide the proof into the following four steps.\\
{\bf Step 1. Group the eigenvalues.} We denote by the eigenvalues of $L_{E,W}^{0}$ by $\left\{e^{-2\pi i\rho_\ell^0}\right\}_{\ell=1}^{2m}$. By Proposition \ref{z}, for any $E\in\mathcal{B}$, there exists $\eta(V)>0$ and $1\leq i_0(E)\leq 2m$, such that we can group $\left\{\rho_\ell^0\right\}_{\ell=1}^{2m}$ as 
$$
E_{1}=\left\{\rho|\rho >\frac{i_0+1}{4m}\eta\right\},\ \ E_{3}=\left\{\rho |\rho<-\frac{i_0+1}{4m}\eta\right\},
$$
$$
E_{2}=\left\{\rho||\rho|\leq \frac{i_0}{4m}\eta\right\},\ \ \#E_{2}\leq 2m_0,
$$
$$
n_1=\#E_{1}=\#E_{3}=n_3.
$$
Thus
$$
\rho>\rho'\geq \frac{\eta}{4m},\ \  \rho\in E_{k}, \rho'\in E_{\ell},\ \ \forall k< \ell.
$$
For any given $\e>0$ and $0<h'<h$, assume $\lambda^{-\frac{1}{2}}\leq \epsilon(\alpha,h,h',m,|L_{E,W}^0|,\e,\eta/4m)$ which is defined in Proposition \ref{final2}. Thus by Proposition \ref{final2}, there exist $B_{j}\in C^{\omega}_{h'}(\T^{d}$,
$GL(m,\C))$, $A_{j}\in GL(m,\C)$ and $f_{j}\in C^{\omega}_{h'}(\T^d,gl(m,\C))$, such that
$$
B_{j}(x+\alpha)L_{E_0,W}^{\lambda^{-1}V}(x)B^{-1}_{j}(x)=A_{j}e^{f_j(x)}.
$$
Moreover, we have the following estimates
\begin{equation}\label{ne1}
|f_{j}|_{h'}\leq \epsilon_{j},\ \ |B_{j}|_0\leq C(\alpha,A)\epsilon_{j-1}^{-200m^2\sigma},
\end{equation}
\begin{equation}\label{ne2}
\left| \Im\rho_{\ell}^j-\Im\rho_{\ell}^{j-1}\right|\leq C(m,A)\epsilon^{\frac{1}{m}-200m^2\sigma}_{j-1},\ \ 1\leq \ell\leq m,
\end{equation}
where $\left\{e^{-2\pi i\rho_\ell^j}\right\}_{\ell=1}^m$ are the eigenvalues of $A_j$ and $\epsilon_j= \epsilon_0^{2^j}$.

By \eqref{ne2}, for $j$ sufficiently large, we can group the eigenvalues of $A_{j+1}$ as $\bigcup\limits_{k=1}^3 E^{j+1}_{k}$ with $\# E^{j+1}_{k}=n_k$, such that
\begin{equation*}
\rho-\rho'\geq \frac{\eta}{8m},\ \ \rho\in E^{j+1}_{k},\rho'\in E^{j+1}_{\ell}, \ \ k< \ell.
\end{equation*}
\noindent{\bf Step 2. Compare $(\alpha, L_{E,W}^{\lambda^{-1}V})$ and $(\alpha, L_{E',W}^{\lambda^{-1}V})$.} Assume  $\delta=|E-E'|$ satisfies $\epsilon_{j+1}\leq \delta\leq \epsilon_{j}$. Then
\begin{align*}
B_{j+1}(x+\alpha)L_{E',W}^{\lambda^{-1} V}B_{j+1}^{-1}(x)&=A_{j+1}e^{f_{j+1}(x)}+B_j(x+\alpha)(L_{E,W}^{\lambda^{-1} V}-L_{E',W}^{\lambda^{-1} V})B_j^{-1}(x)\\
&:=A_{j+1}e^{\tilde{f}_{j+1}}.
\end{align*}
By \eqref{ne1}, we have
$$
|\tilde{f}_{j+1}|_h\leq 100\delta|B_{j+1}|_0^2+\epsilon_{j+1}\leq 100\delta\epsilon_j^{-800m^2\e}\leq 100\delta^{1-800m^2\e}.
$$ 
\noindent{\bf Step 3. H\"older continuity of the Lyapunov exponents.}
Recall that $\{L_\ell(E)\}_{\ell=1}^m$ are the nonnegative Lyapunov exponents of $(\alpha, L_{E,W}^{\lambda^{-1}V})$. By Corollary \ref{final3}, we have
\begin{equation*}
\left|\sum\limits_{\ell=1}^{n_1}L_\ell(E)-\sum\limits_{\ell=1}^{n_1}\Im\rho^{j+1}_{\ell}\right|\leq C(\eta,V)\delta^{\frac{1}{2}},
\end{equation*}
\begin{equation*}
\left|L_\ell(E)-\Im\rho^{j+1}_{\ell}\right|\leq C(\eta,V)\delta^{\frac{1}{n_2}-800m^2\e}, \ \  n_1+1\leq \ell\leq n_1+n_2,
\end{equation*}
\begin{equation*}
\left|\sum\limits_{\ell=1}^{n_1}L_\ell(E')-\sum\limits_{\ell=1}^{n_1}\Im\rho^{j+1}_{\ell}\right|\leq C(\eta,V)\delta^{\frac{1}{2}},
\end{equation*}
\begin{equation*}
\left|L_\ell(E')-\Im\rho^{j+1}_{\ell}\right|\leq C(\eta,V)\delta^{\frac{1}{n_2}-800m^2\e}, \ \  n_1+1\leq \ell\leq n_1+n_2.
\end{equation*}
It follows
\begin{equation}\label{neq1}
\left|\sum\limits_{\ell=1}^{n_1}L_\ell(E)-\sum\limits_{\ell=1}^{n_1}L_\ell(E')\right|\leq C(\eta,V)\delta^{\frac{1}{2}},
\end{equation}
\begin{equation}\label{neq2}
\left|L_\ell(E')-L_\ell(E)\right|\leq C(\eta,V)\delta^{\frac{1}{n_2}-800m^2\e}, \ \  n_1+1\leq \ell\leq n_1+n_2.
\end{equation}
\eqref{neq1} and \eqref{neq2} imply
\begin{align*}
\left|\sum\limits_{\ell=1}^m L_\ell(E)-\sum\limits_{\ell=1}^m L_\ell(E')\right|&\leq \left|\sum\limits_{\ell=1}^{n_1} L_\ell(E)-\sum\limits_{\ell=1}^{n_1}L_\ell(E')\right|+\left|\sum\limits_{\ell=n_1+1}^m L_\ell(E)-\sum\limits_{\ell=n_1+1}^m L_\ell(E')\right|\\
&\leq C(\eta,V)\delta^{\frac{1}{2}}+2mC(\eta,V)\delta^{\frac{1}{2m_0}-800m^2\e}\\
&\leq  4mC(\eta,V)|E-E'|^{\frac{1}{2m_0}-800m^2\e}.
\end{align*}
\noindent {\bf Step 4. H\"older continuity of the  IDS.} By Thouless formula 
$$
\sum\limits_{\ell=1}^m L_\ell(E)+\ln|W_m|=\int \ln\lvert E-E'\rvert dN(E'),
$$
we have
\begin{flalign*}
L(E+i\epsilon)-L(E)=&\frac{1}{2}\int \ln (1+\frac{\epsilon^2}{(E-E^{'})^2})dN(E')\\
\geq &c'|N(E+\epsilon)-N(E-\epsilon)|.
\end{flalign*}
Thus
$$
|N(E+\epsilon)-N(E-\epsilon)|\leq C(\alpha,V,W)\epsilon^{\frac{1}{2m_0}-800m^2\e}.
$$
\section{Appendix}
In the appendix, we give the some useful lemmas.
\begin{Lemma}\label{useful_lem}
Assume 
$$
A=\begin{pmatrix}a_{11}&\cdots&a_{1n}\\
&\ddots&\vdots\\ && a_{nn}\end{pmatrix},\  \ 
B=\begin{pmatrix}b_{11}&\cdots&b_{1m}\\
&\ddots&\vdots\\ && b_{mm}\end{pmatrix},
$$ 
\begin{align}\label{noncond}
|a_{ii}e^{\pm2\pi i\langle k,\alpha\rangle}-b_{jj}| \geq \eta, \ \ 1\leq i\leq n, \ \ 1\leq j\leq m,\ \  \forall k\in\Lambda.
\end{align}
Then for any $Y\in\{f|f(x)=\sum\limits_{k\in\Lambda}\hat{f}(k)e^{2\pi i\langle k,x\rangle}\}$, we have
$$
\lvert Y(\cdot+\alpha)A-BY(\cdot)\rvert_h\geq \eta\left(1+\max\{m,n\}(|A|+|B)|\eta^{-1}\right)^{-(m+n)}\lvert Y\rvert_h.
$$
\end{Lemma}
\begin{pf}
Without loss of generality, we assume that $n\leq m$. Let $F(x)=Y(x+\alpha)A-BY(x)$, we inductively prove for $k\in\Lambda$,
\begin{equation}\label{f1}
|\widehat{Y}_{i,1}(k)|\leq \eta^{-1}\left(1+m(|A|+|B|)\eta^{-1}\right)^{m-i}|\widehat{F}(k)|.
\end{equation}
\begin{equation}\label{f3}
|\widehat{Y}_{m,j}(k)|\leq \eta^{-1}\left(1+m(|A|+|B|)\eta^{-1}\right)^{j-1}|\widehat{F}(k)|.
\end{equation}       
\begin{equation}\label{f2}
|\widehat{Y}_{i,j}(k)|\leq \eta^{-1}\left(1+m(|A|+|B)|\eta^{-1}\right)^{m-i+j-1}|\widehat{F}(k)|.
\end{equation}
We first prove \eqref{f1}, for $i=m$, we have
$$
F_{m,1}(x)=a_{11}Y_{m,1}(x+\alpha)-b_{mm}Y_{m,1}(x),
$$
thus
$$
\widehat{F}_{m,1}(k)=a_{11}e^{2\pi i\langle k,\alpha\rangle}\widehat{Y}_{m,1}(k)-b_{mm}\widehat{Y}_{m,1}(k),
$$
by \eqref{noncond}, we have
$$
|\widehat{Y}_{m,1}(k)|\leq \eta^{-1}|\widehat{F}_{m,1}(k)|\leq \eta^{-1}|\widehat{F}(k)|.
$$
Assume for $i_0\leq i\leq m$, we already have
\begin{equation*}
|\widehat{Y}_{i,1}(k)|\leq \eta^{-1}\left(1+m(|A|+|B|)\eta^{-1}\right)^{m-i}|\widehat{F}(k)|.
\end{equation*}
Then
$$
F_{i_0-1,1}(x)=a_{11}Y_{i_0-1,1}(x+\alpha)-b_{i_0-1,i_0-1}Y_{i_0-1,1}(x)+\sum\limits_{\ell=i_0}^{m}b_{i_0-1,\ell}Y_{\ell,1}(x).
$$
By \eqref{noncond}, we have
$$
|\widehat{F}_{i_0-1,1}(k)|\geq \eta^{-1} |\hat{Y}_{i_0-1,1}(k)|-|B|\sum\limits_{\ell=i_0}^{m}|\widehat{Y}_{\ell,1}(k)|,
$$
this implies that
\begin{align*}
|\widehat{Y}_{i_0-1,1}(k)|&\leq \eta^{-1}\left(|\widehat{F}_{i_0-1,1}(k)|+|B|\sum\limits_{\ell=i_0}^{m}|\widehat{Y}_{\ell,1}(k)|\right)\\
&\leq \eta^{-1}\left(1+|B|\sum\limits_{\ell=i_0}^{m}\eta^{-1}(1+m(|A|+|B|)\eta^{-1})^{m-\ell}\right)|\hat{F}(k)|\\
&\leq \eta^{-1}\left(1+m|B|\eta^{-1}(1+m(|A|+|B|)\eta^{-1})^{m-i_0}\right)|\widehat{F}(k)|\\
&\leq \eta^{-1}\left(1+m(|A|+|B|)\eta^{-1}\right)^{m-i_0+1}|\widehat{F}(k)|.
\end{align*}
The proof of \eqref{f3} is exactly the same.

Now, we inductively prove \eqref{f2}, assume for $(i,j)\in \mathcal{B}_{i_0,j_0}$ where
$$
\mathcal{B}_{i_0,j_0}=\{(i,j)|i_0 \leq i\leq n,1\leq j\leq m\}\cup\{(i,j)|1 \leq i\leq m,1\leq j\leq j_0\},
$$
we already have
\begin{equation*}
|\widehat{Y}_{i,j}(k)|\leq \eta^{-1}\left(1+m(|A|+|B|)\eta^{-1}\right)^{m-i+j-1}|\widehat{F}(k)|.
\end{equation*}
Then
\begin{align*}
F_{i_0-1,j_0+1}(x)=&-b_{i_0-1,i_0-1}Y_{i_0-1,j_0+1}(x)+a_{j_0+1,j_0+1}Y_{i_0-1,j_0+1}(x+\alpha)\\
&-\sum\limits_{\ell=i_0}^{m}b_{i_0-1,\ell}Y_{\ell,j_0+1}(x)+\sum\limits_{\ell=1}^{j_0}a_{\ell,j_0+1}Y_{i_0-1,\ell}(x+\alpha).
\end{align*}
By \eqref{noncond}, we have
$$
|\widehat{F}_{i_0-1,j_0+1}(k)|\geq \eta |\widehat{Y}_{i_0-1,j_0+1}(k)|-|B|\sum\limits_{\ell=i_0}^{m}|\widehat{Y}_{\ell,j_0+1}(k)|-|A|\sum\limits_{\ell=1}^{j_0}|\widehat{Y}_{i_0-1,\ell}(k)|,
$$
this implies that
\begin{align*}
|\widehat{Y}_{i_0-1,j_0+1}(k)|&\leq \eta^{-1}\left(|\widehat{F}_{i_0-1,j_0+1}(k)|+|B|\sum\limits_{\ell=i_0}^{m}|\widehat{Y}_{\ell,j_0+1}(k)|+|A|\sum\limits_{\ell=1}^{j_0}|\hat{Y}_{i_0-1,\ell}(k)|\right)\\
&\leq \eta^{-1}\left(1+|B|\sum\limits_{\ell=i_0}^{m}\eta^{-1}(1+m(|A|+|B|)\eta^{-1})^{m-\ell+j_0}\right)|F(k)|\\
&\ \ \ \ +\left(|A|\sum\limits_{\ell=1}^{j_0}\eta^{-1}(1+m(|A|+|B|)\eta^{-1})^{m-i_0+\ell}\right)|F(k)|\\
&\leq \eta^{-1}\left(1+m(|A|+|B|)\eta^{-1}(1+m(|A|+|B|)\eta^{-1})^{m-i_0+j_0}\right)|F(k)|\\
&\leq \eta^{-1}\left(1+m(|A|+|B|)\eta^{-1}\right)^{m-i_0+1+j_0}|F(k)|.
\end{align*}
Similarly, we can prove that for any $(i,j)\in \mathcal{B}_{i_0-1,j_0+1}$,
\begin{equation*}
|\widehat{Y}_{i,j}(k)|\leq \eta^{-1}\left(1+m(|A|+|B|)\eta^{-1}\right)^{m-i+j-1}|\hat{F}(k)|.
\end{equation*}
Thus we finish the proof of \eqref{f2}. By the definition of analytic norm, we have
$$
\lvert Y(\cdot+\alpha)A-BY(\cdot)\rvert_h\geq \eta^{-1}\left(1+m(|A|+|B)|\eta^{-1}\right)^{-(m+n)}\lvert Y\rvert_h.
$$
\end{pf}
For any given $\eta>0$, $\alpha\in \R^{d}$ and $A\in GL(m,\C)$, we decompose $\mathcal{B}_h=C^{\omega}_{h}(\T^{d},gl(m,\C))=\mathcal{B}_h^{nre}(\eta) \bigoplus\mathcal{B}_h^{re}(\eta)$ in such a way that for any $Y\in\mathcal{B}_h^{nre}(\eta)$,
\begin{equation}\label{space}
A^{-1}Y(\theta+\alpha)A\in\mathcal{B}_h^{nre}(\eta), \qquad \lvert A^{-1}Y(\theta+\alpha)A-Y(\theta)\rvert_h\geq\eta\lvert Y(\theta)\rvert_h.
\end{equation}
Moreover, let $\mathbb{P}_{nre}$ and $\mathbb{P}_{re}$ be the standard projections from $\mathcal{B}_h$ onto $\mathcal{B}_h^{nre}(\eta)$ and $\mathcal{B}_h^{re}(\eta)$ respectively. 
\begin{Lemma}[Lemma 3.1 of \cite{ccyz}]\label{lem2} Assume that $\epsilon\leq (4 |A|)^{-4}$ and  $\eta \geq 13| A|^2{\epsilon}^{\frac{1}{2}}$. For any $g\in \mathcal{B}_h$ with $|g|_h \leq \epsilon$,  there exist $Y\in \mathcal{B}_h$ and $g^{re}\in \mathcal{B}_h^{re}(\eta)$ such that
$$
e^{Y(x+\alpha)}(Ae^{g(x)})e^{-Y(x)}=Ae^{g^{re}(x)},
$$
with $\lvert Y \rvert_h\leq \epsilon^{\frac{1}{2}}$ and $\lvert g^{re}\rvert_h\leq 2\epsilon$.
\end{Lemma}
\begin{Remark}
Although in Lemma 3.1 in \cite{ccyz}, the authors assume that $A\in SU(1,1)$,  the proof only uses implicit function theorem and essentially works for any $A\in GL(m,\C)$. See the continuous version of Lemma \ref{lem2} in \cite{hy}.
\end{Remark}

\begin{Proposition}\label{z}
Assume the zeros (counting multiplicity) of $W(\theta)-E$ on $\T$ are no more than $2m_0$ for all $E\in \R$, let $\mathcal{Z}(E)=\{z\in\mathbb{C}| W(z)-E=0\}$. There exist $\eta(V)>0$ and $1\leq i_0(E)\leq 2m$, such that $\mathcal{Z}(E)=\mathcal{Z}^+(E)\cup\mathcal{Z}^-(E)\cup \mathcal{Z}_0(E)$ where
$$
\mathcal{Z}^+(E)=\left\{z\in\mathcal{Z}(E)|\ln|z|>\frac{i_0+1}{4m}\eta\right\},\ \ \mathcal{Z}^-(E)=\left\{z\in\mathcal{Z}(E)|\ln|z|<-\frac{i_0+1}{4m}\eta\right\},
$$
$$
 \mathcal{Z}_0(E)=\left\{z\in\mathcal{Z}(E)||\ln|z||\leq \frac{i_0}{4m}\eta\right\},\ \ |\mathcal{Z}_0(E)|\leq 2m_0,
$$
$$
|\mathcal{Z}^+(E)|=|\mathcal{Z}^-(E)|.
$$
\end{Proposition}
\begin{pf}
Let $\mathcal{Z}(E)=\{z_i(E)\}_{i=1}^{2m}$ satisfying $|z_1(E)|\geq |z_2(E)|\geq \cdots\geq |z_{2m}(E)|$. Since the zeros (counting multiplicity) of $W(\theta)-E$ on $\T$ are no more than $2m_0$ for all $E\in \R$, we have $|z_{m-m_0}(E)|>1$ for all $E\in \R$. On the other hand, the zeros the polynomial depend continuously on $E$. Thus there exists $\eta(V)>0$, such that $|z_{m-m_0}(E)|\geq 1+\eta$ for all $E\in\R$. Now we fix $E$, by pigeonhole principle, there exists $1\leq i_0\leq 2m$ such that $\{|\ln|z||: z\in\mathcal{Z}_0\}\cap [\frac{i_0}{3m}\eta,\frac{i_0+1}{3m}\eta]=\emptyset$. We denote by 
$$
\mathcal{Z}^+(E)=\left\{z\in\mathcal{Z}(E)|\ln|z|>\frac{i_0+1}{3m}\eta\right\},\ \ \mathcal{Z}^-(E)=\left\{z\in\mathcal{Z}(E)|\ln|z|<-\frac{i_0+1}{3m}\eta\right\},
$$
$$
 \mathcal{Z}_0(E)=\left\{z\in\mathcal{Z}(E)||\ln|z||\leq \frac{i_0}{3m}\eta\right\},\ \ |\mathcal{Z}_0(E)|\leq 2m_0,
$$
then $\mathcal{Z}(E)=\mathcal{Z}^+(E)\cup\mathcal{Z}^-(E)\cup \mathcal{Z}_0(E)$.
\end{pf}

\section*{Acknowledgement}
 L. Ge was partially supported by NSF DMS-190146 and AMS-Simons Travel Grant 2020–2022. J. You was partially supported by NNSF of China (11871286) and
Nankai Zhide Foundation. X. Zhao was partially supported by NSF DMS-190146 and China Scholarship Council (No. 201906190072).

\end{document}